\documentclass[aap,seceqn,MSNbibl,citesort,dvips]{arximspdf}
\usepackage[vtexbibl]{}

\doi{10.1214/09-AAP627}
\volume{20}
\issue{2}
\pubyear{2010}
\firstpage{565}
\lastpage{592}

\makeatletter

\newtheorem{propo}{Proposition}[section]
\newtheorem{lemma}[propo]{Lemma}
\newproclaim{definition}[propo]{Definition}
\newtheorem{coro}[propo]{Corollary}
\newtheorem{theorem}[propo]{Theorem}
\newproclaim{Remark}{Remark}
\newproclaim{remark}[propo]{Remark}

\makeatother

\begin{document}
\begin{frontmatter}

\title{Ising models on locally tree-like graphs\thanksref{T1}}
\runtitle{Ising models on locally tree-like graphs}
\thankstext{T1}{Supported in
part by NSF Grant DMS-08-06211.}
\begin{aug}
\author[A]{\fnms{Amir} \snm{Dembo}}
\and
\author[B]{\fnms{Andrea} \snm{Montanari}\ead[label=e2]{montanari@stanford.edu}\corref{}}
\runauthor{A. Dembo and A. Montanari}
\affiliation{Stanford University}
\address[A]{Departments of Statistics\\
\quad and Mathematics\\
Stanford University\\
Stanford, California 94305-9505\\
USA} 
\address[B]{Departments of Electrical\\\quad
Engineering and Statistics\\ Stanford University\\
Stanford, California 94305-9505\\
USA\\
\printead{e2}}
\end{aug}

\received{\smonth{5} \syear{2008}}
\revised{\smonth{7} \syear{2009}}

%
\begin{abstract}
We consider ferromagnetic Ising models on graphs
that converge locally to trees. Examples include random regular
graphs with bounded degree and uniformly random graphs with
bounded average degree. We prove that the ``cavity''
prediction for the
limiting
free energy per spin is correct for
\emph{any positive temperature and external field}.
Further, local marginals can be approximated
by iterating a set of mean field (cavity) equations.
Both results are achieved by proving the local convergence
of the Boltzmann distribution on the original graph to the Boltzmann
distribution on the appropriate infinite random tree.
\end{abstract}

%

\end{frontmatter}

\setcounter{footnote}{1}
\section{Introduction}

A ferromagnetic \emph{Ising model on the
finite
graph $G$} (with vertex set $V$,
and edge set $E$) is defined by the following Boltzmann distributions over
$\underline{x}= \{x_i\dvtx  i\in V\}$,
with
$x_i\in\{+1,-1\}$:
%
\begin{eqnarray}\label{eq:IsingModel}
\mu(\underline{x}) =\frac{1}{Z(\beta,B)}  \exp\biggl\{\beta\sum
_{(i,j)\in E}x_ix_j
+B\sum_{i\in V}x_i\biggr\}  .
\end{eqnarray}
These distributions are parametrized by the
``magnetic field'' $B$ and ``inverse temperature'' $\beta\ge0$,
where the partition function $Z(\beta,B)$ is
fixed by the normalization condition $\sum_{\underline{x}} \mu
(\underline{x})=1$.
Throughout the paper, we will be interested in sequences of
graphs\footnote{We adopt the notation
$[i]=\{1,2,\dots,i\}$ for the set of first $i$ integers.}
$G_n = (V_n\equiv[n],E_n)$ of diverging size $n$.

Nonrigorous statistical mechanics techniques, such as the ``replica''
and ``cavity methods,'' allow to make a number of predictions
on the model (\ref{eq:IsingModel}), when the graph $G$
``lacks
any finite-dimensional structure.''
The most basic quantity in this context is the
asymptotic \emph{free entropy density}
%
\begin{eqnarray}\label{eq:FreeEnergy}
\phi(\beta,B) \equiv\lim_{n\to\infty}\frac{1}{n}\log Z_n(\beta
,B)
\end{eqnarray}
(this quantity is also
sometimes called in the literature also free energy or pressure).
The limit free entropy density and the large deviation properties
of Boltzmann distribution were characterized
in great detail \cite{NewmanEllis} in the case of a complete
graph $G_n=K_n$ (the inverse temperature must then be scaled by $1/n$
to get a nontrivial limit).
Statistical physics predictions exist, however, for a much wider class
of graphs,
including most notably sparse random graphs with bounded average degree;
see, for instance, \cite{Dorogotsev,Johnston,Vespignani}.
This is a direction of interest for at least two reasons:
\begin{longlist}[(ii)]
\item[(i)]
Sparse graphical structures arise in a number of problems from
combinatorics and theoretical computer science.
Examples include random satisfiability, coloring of random graphs and
graph partitioning \cite{MezardMontanari}.
In all of these cases, the uniform measure
over solutions can be regarded as the Boltzmann distribution
for a modified spin glass with multispin interactions.
Such problems have been successfully attacked using nonrigorous
statistical mechanics techniques.

A mathematical foundation of this approach is still lacking, and would be
extremely useful.

\item[(ii)]
Sparse graphs allow to introduce a nontrivial notion of distance
between vertices, namely the length of the shortest path
connecting them. This geometrical structure allows for new characterizations
of the measure (\ref{eq:IsingModel}) in terms of correlation decay.
This type of characterization is in turn related to the theory of Gibbs
measures on infinite trees \cite{KMRSZ}.
\end{longlist}

The asymptotic free entropy density (\ref{eq:FreeEnergy}) was
determined rigorously only in a few cases for sparse graphs.
In \cite{GerschenMon1}, this task was accomplished for random regular graphs.
De Sanctis and Guerra \cite{deSantisGuerra} developed interpolation
techniques for random graphs with independent edges (Erd\"{o}s--Renyi type)
but only determined the free entropy density at high temperature and at zero
temperature (in both cases with vanishing magnetic field).
The latter is in fact equivalent to counting the number of
connected components of a random graph.
Interestingly, the partition function $Z_n(\beta,B)$
can be approximated in polynomial time for $\beta\ge0$, using an appropriate
Markov chain Monte Carlo algorithm \cite{JerrumSinclair}.
It is intriguing that no general approximation algorithms exists
in the case $\beta<0$ (the ``antiferromagnetic'' Ising model).
Correspondingly, the statistical physics conjecture for the
free entropy density \cite{MezardMontanari} becomes significantly
more intricate (presenting the so-called ``replica symmetry breaking''
phenomenon).

In this paper we generalize the previous results by
rigorously verifying the validity of the
Bethe free entropy prediction for the value of the
limit in (\ref{eq:FreeEnergy}) for generic graph
sequences that converge locally to trees.
Indeed, we control the
free entropy density by proving that the Boltzmann measure
(\ref{eq:IsingModel})
converges locally to the Boltzmann measure of a model on
a tree. The philosophy is related
to the local weak convergence method of \cite{Aldous}.

Finally, several of the proofs have an algorithmic interpretation,
providing an efficient procedure for approximating the local marginals
of the
Boltzmann measure. The essence of this procedure consists in solving
by iteration certain mean field (cavity) equations.
Such an algorithm is known in artificial intelligence and computer
science under the name of \emph{belief propagation}. Despite its success
and wide applicability, only weak performance guarantees
have been proved so far. Typically, it is possible to prove its correctness
in the high temperature regime, as a consequence
of a uniform decay of correlations holding there (spatial mixing)
\cite{CompTree,Gamarnik,Devavrat}.
The behavior of iterative inference algorithms on Ising models
was recently considered in \cite{Mooij,Sudderth}.

The emphasis of the present paper is on
the low-temperature regime in which uniform decorrelation does not hold.
We are able to prove that belief propagation converges exponentially
fast on any graph, and that the resulting estimates are asymptotically
exact for large locally tree-like graphs. The main idea is to
introduce a magnetic field to break explicitly the $+/-$ symmetry,
and to carefully exploit the monotonicity properties of the model.

A key step consists of estimating the correlation between the root
spin of an Ising model on a tree and positive boudary conditions.
Ising models on trees are interesting per se, and
have been the object of significant mathematical work; see, for instance,
\cite{Lyons,Steif,PeresReconstr}. The question considered here
appears, however, to be novel.

The next section provides the basic technical definitions (in particular
concerning graphs and local convergence to trees), and the formal
statement of our main results.
Notation and certain key tools are described in Section \ref{sec:Notations}
with Section \ref{sec:TreeProofs}
devoted to proofs of the relevant properties of Ising models on trees
(which are of independent interest). The latter are used in
Sections \ref{sec:Algo} and \ref{sec:Graphs}
to derive our main results concerning models on tree-like graphs.
A companion paper \cite{GerschenMon2} deals with the related
challenging problem of spin glass models on sparse graphs.
%
%
\section{Definitions and main results}
\label{sec:DefinitionMain}

The next subsections contain some basic definitions on graph
sequences and the notion of local convergence to random trees.
Sections \ref{sec:FreeEntropy} and \ref{sec:AlgoImpl} present
our results on the free entropy density and the algorithmic implications
of our analysis.
%
%
\subsection{Locally tree-like graphs}
\label{sec:Local}

Let $P= \{P_k\dvtx   k\ge0\}$ a probability distribution over
the nonnegative integers, with finite,
positive first moment, and denote by
%
\begin{eqnarray}
\rho _k = \frac{kP_k}{\sum_{l=1}^{\infty}lP_l}  ,
\end{eqnarray}
its size-biased version. For any $t\ge0$, we let
$\mathsf{T}(P,\rho ,t)$ denote the random rooted tree generated as follows.
First draw an integer $k$ with distribution $P_k$, and connect
the root to $k$ offspring. Then recursively, for each node in
the last generation, generate an integer $k$ independently with distribution
$\rho _k$, and connect the node to $k-1$ new nodes. This is repeated
until the tree has $t$ generations.

Sometimes it will be useful to consider the ensemble $\mathsf{T}(\rho ,t)$
whereby the root node has degree $k-1$ with probability $\rho _k$.
We will drop the degree distribution arguments from
$\mathsf{T}(P,\rho ,t)$ or $\mathsf{T}(\rho ,t)$ and write $\mathsf
{T}(t)$ whenever
clear from the context.
Notice that the infinite trees
$\mathsf{T}(P,\rho ,\infty)$ and $\mathsf{T}(\rho ,\infty)$ are
well defined.

The average branching factor of trees will be denoted by $\overline
{\rho}$,
and the average root degree by $\overline{P}$. In formulae
%
\begin{eqnarray}
\overline{P}\equiv\sum_{k=0}^{\infty}k P_k   ,\qquad
\overline{\rho}\equiv\sum_{k=1}^{\infty}(k-1)   \rho _k   .
\end{eqnarray}

We denote by $G_n=(V_n,E_n)$ a graph with vertex set
$V_n\equiv[n]=\{1,\dots,n\}$. The distance $d(i,j)$ between
$i,j\in V_n$ is the length of the shortest path from $i$ to $j$ in~$G_n$. Given a vertex $i\in V_n$, we let
$\mathsf{B}_i(t)$ be the set of vertices whose distance from~$i$ is at most
$t$. With a slight abuse of notation, $\mathsf{B}_i(t)$ will also
denote the
subgraph induced by those vertices. For $i\in V_n$, we let ${\partial
i}$ denote
the set of its neighbors ${\partial i}\equiv\{j\in V_n\dvtx  (i,j)\in
E_n\}$,
and $|{\partial i}|$ its size (i.e. the degree of $i$).

This paper is concerned by sequence of graphs $\{G_n\}_{n\in\mathbb{N}}$
of diverging size, that converge locally to trees.
Consider two trees $T_1$ and $T_2$ with vertices labeled arbitrarily.
We shall write $T_1\simeq T_2$ if the two trees become identical when vertices
are relabeled from $1$ to $|T_1|=|T_2|$, in a
breadth first fashion, and following lexicographic order among siblings.
\begin{definition}
Considering a sequence of graphs $\{G_n\}_{n\in\mathbb{N}}$, let
$\mathbb{P}_n$ denote the law induced on the ball $\mathsf{B}_i(t)$
in $G_n$
centered at a uniformly chosen random vertex $i\in[n]$.
We say that $\{G_n\}$ \emph{converges locally} to the random tree
$\mathsf{T}(P,\rho ,\infty)$ if, for any $t$, and any rooted tree $T$
with $t$ generations
%
\begin{eqnarray}
\lim_{n\to\infty}\mathbb{P}_n\{\mathsf{B}_i(t) \simeq T\} =
\mathbb{P}\{
\mathsf{T}(P,\rho,t) \simeq T\}  .
\end{eqnarray}
\end{definition}
\begin{definition}
We say that a sequence of graphs
$\{G_n\}_{n\in\mathbb{N}}$ is \emph{uniformly sparse} if
%
\begin{eqnarray}
\lim_{l \to\infty} \limsup_{n \to\infty}
\frac{1}{n} \sum_{i\in V_n}|{\partial i}|  \mathbb{I}(|{\partial
i}|\ge l) = 0   .
\end{eqnarray}
\end{definition}
%

\subsection{Free entropy}
\label{sec:FreeEntropy}

According to the statistical physics derivation \cite{Vespignani}, the model
(\ref{eq:IsingModel})
has a line of first-order phase
transitions for $B=0$ and $\beta>\beta_{\rm c}$
[i.e., where the continuous function $B \mapsto\phi(\beta,B)$
exhibits a discontinuous derivative].
The critical temperature depends on the
graph only through the average branching factor and is determined by the
condition
%
\begin{eqnarray}
\overline{\rho} \tanh\beta_{\rm c} = 1  .
\end{eqnarray}
Notice that $\beta_{\rm c}\simeq1/\overline{\rho}$ for large degrees.

The asymptotic free-entropy density is given in terms of the
fixed point of a distributional recursion.
One characterization of this fixed point is as follows.
\begin{lemma}\label{lemma:Recursive}
Consider the sequence of random variables $\{h^{(t)}\}$ defined by $h^{(0)}=0$
identically and, for $t\ge0$,
%
\begin{equation}\label{eqn:h_recursion}
h^{(t+1)} \stackrel{\mathrm{d}}{=}B + \sum_{i=1}^{K-1} \xi\bigl(\beta
,h^{(t)}_i\bigr)   ,
\end{equation}
where $K$ is an integer valued random variable of distribution $\rho $,
%
\begin{equation}\label{eq:xidef}
\xi(\beta,h) \equiv\operatorname{atanh}[\tanh(\beta) \tanh(h)] ,
\end{equation}
and the $h^{(t)}_i$\textup{'}s
are i.i.d. copies of $h^{(t)}$ that are independent of $K$.
If $B>0$ and $\rho $ has finite first moment, then the
distributions of $h^{(t)}$ are stochastically monotone
and $h^{(t)}$ converges in distribution to the unique fixed point $h^*$
of the recursion (\ref{eqn:h_recursion}) that is
supported on $[0,\infty)$.
\end{lemma}

Our next result confirms the statistical physics prediction for the
free-entropy density.
\begin{theorem}\label{theorem:free_energy}
Let $\{G_n\}_{n\in\mathbb{N}}$ be a sequence of uniformly sparse graphs
that converges locally to $\mathsf{T}(P,\rho ,\infty)$. If $\rho $
has finite first moment (that is if $P$ has finite second moment),
then for any $B \in{\mathbb{R}}$ and $\beta\ge0$ the following
limit exists:
%
\begin{equation}\label{eqn:free_energy}
\lim_{n\to\infty}\frac{1}{n}\log Z_n(\beta,B) =\phi(\beta,B)  .
\end{equation}
Moreover, for $B>0$ the limit is given by
%
\begin{eqnarray}\label{eqn:phi}
\phi(\beta,B) &\equiv& \frac{\overline{P}}{2}  \log\cosh(\beta
) -
\frac{\overline{P}}{2}  \mathbb{E}\log[1+\tanh(\beta) \tanh
(h_1)\tanh(h_2)]
\nonumber
\\
&&{}+ \mathbb{E} \log\Biggl\{ e^B \prod_{i=1}^{L} [1 +
\tanh(\beta) \tanh(h_i)]
\\
&&{}\qquad\hspace{18pt} + e^{-B} \prod_{i=1}^{L}
[1 - \tanh(\beta) \tanh(h_i) ] \Biggr\}  ,\nonumber
\end{eqnarray}
where $L$ has distribution $P_{l}$ and is independent of the
\textup{``}cavity fields\textup{''} $h_i$ that are i.i.d. copies of the fixed point
$h^*$ of Lemma \ref{lemma:Recursive}.
Also, $\phi(\beta,B)=\phi(\beta,-B)$
and $\phi(\beta,0)$ is the limit of $\phi(\beta,B)$ as $B \to0$.
\end{theorem}

The proof of Theorem \ref{theorem:free_energy} is based on two steps:
\begin{longlist}[(a)]
\item[(a)] Reduce the computation of
$\phi_n(\beta,B) =\frac{1}{n} \log Z_n(\beta,B)$ to computing expectations
of local (in $G_n$)
quantities with respect to the Boltzmann measure (\ref{eq:IsingModel}).
This is achieved by noticing that the derivative of $\phi_n(\beta,B)$
with respect to $\beta$ is a sum of such expectations.
\item[(b)] Show that expectations of local quantities on $G_n$ are well
approximated by the same expectations with respect to an Ising model on the
associated tree $\mathsf{T}(P,\rho ,t)$ (for $t$ and $n$ large).
This is proved by showing that, on such a tree, local expectations are
insensitive to boundary conditions that dominate stochastically
free boundaries.
The theorem then follows by monotonicity arguments.
\end{longlist}
The key step is of course the last one. A stronger requirement would be that
these
expectation values are insensitive to any boundary condition,
which would coincide with uniqueness of the Gibbs measure on
$\mathsf{T}(P,\rho ,\infty)$. Such a requirement would allow
for an elementary proof, but holds only at ``high'' temperature,
$\beta\le\beta_{\rm c}$.

Indeed, insensitivity to positive boundary conditions is proved
in Section \ref{sec:TreeProofs} for
the following collection of
trees of conditionally independent (and of bounded
average) offspring numbers.
\begin{definition}\label{def-cit}
An infinite tree $\mathsf{T}$ rooted at the vertex $\o$
is called \emph{conditionally independent} if
for each integer $k \ge0$, conditional on
the subtree $\mathsf{T}(k)$ of the first $k$ generations of $\mathsf{T}$,
the number of offspring $\Delta_j$ for $j \in\partial\mathsf{T}(k)$
are independent of each other, where
$\partial\mathsf{T}(k)$ denotes the set of vertices at generation $k$.
We further assume that the [conditional on $\mathsf{T}(k)$]
first moments of $\Delta_j$ are uniformly bounded by a given
nonrandom finite constant $\Delta$.
\end{definition}

Beyond the random tree $\mathsf{T}(P,\rho ,\infty)$,
these include deterministic trees with bounded degrees
and certain multi-type branching processes (such as
random bipartite trees and percolation clusters
on deterministic trees of bounded degree).
Consequently, Theorem \ref{theorem:free_energy}
extends to any uniformly sparse graph sequence
that converge locally to a random tree $\mathsf{T}$ of
the form of Definition \ref{def-cit} except that
the formula $\phi(\beta,B)$ is in general more
involved than the one given in (\ref{eqn:phi}).
For example, such an extension allows one to handle
uniformly random bipartite graphs with different
degree distributions $P_k$ and $Q_k$ for the two types of vertices.

While we refrain from formalizing and proving such
generalizations, we note in passing that our derivation
of the formula (\ref{eqn:phi}) implicitly
uses the fact that $\mathsf{T}(P,\rho ,\infty)$ possesses the
involution invariance of \cite{Aldous}. As pointed out
in \cite{AL}, every local limit of finite graphs
must have the involution invariance property (which
clearly not every conditionally independent tree has).

%
%
\subsection{Algorithmic implications}
\label{sec:AlgoImpl}

The free entropy density is not the only quantity that can be
characterized for Ising models on locally tree-like graphs.
Indeed local marginals can be efficiently computed with good accuracy.
The basic idea is to solve a set of mean field equations iteratively.
These are known as Bethe--Peierls or cavity equations and the
corresponding algorithm is referred to as ``belief propagation'' (BP).

More precisely, associate to each directed edge in the graph $i\to j$,
with $(i,j)\in G$, a distribution $\nu_{i\to j}(x_i)$ over $x_i\in
\{+ 1,-1\}$.
In the computer science literature these distributions are referred
to as ``messages.'' They are updated as follows:
%
\begin{eqnarray}\label{eq:BPIteration}
\nu_{i\to j}^{(t+1)}(x_i) =\frac{1}{z^{(t)}_{i\to j}}  e^{Bx_i}
\prod_{l\in{\partial i}\setminus j}\sum_{x_l}e^{\beta x_ix_l}\nu
_{l\to i}^{(t)}
(x_l) .
\end{eqnarray}
The initial conditions $\nu_{i\to j}^{(0)}( \cdot )$ may be taken
to be
uniform or chosen\vspace{-2pt} according to some heuristic. We will say that the
initial condition is \emph{positive}\vspace{-2pt} if $\nu_{i\to j}^{(0)}(+1)\ge
\nu_{i\to j}^{(0)}(-1)$ for each of these messages.\vspace{1pt}

Our next result concerns the uniform exponential
convergence of the BP iteration to the same fixed
point of (\ref{eq:BPIteration}), irrespective of its
positive initial condition.
Here and below, we denote by $\Vert p-q\Vert_{\mathrm{TV}}$ the total
variation distance
between distributions $p$ and $q$.
\begin{theorem}\label{theorem:BPConvergence}
Assume $\beta\ge0$, $B>0$ and $G$ is a graph of finite
maximal degree~$\Delta$. Then,
there exists $A=A(\beta,B,\Delta)$ finite,
$\lambda=\lambda(\beta,B,\Delta)>0$
and a fixed point $\{\nu^*_{i\to j}\}$ of the BP iteration
(\ref{eq:BPIteration}) such that for any positive initial\vspace{-2pt} condition
$\{\nu^{(0)}_{l \to k}\}$ and all $t \ge0$,
%
\begin{eqnarray}
\label{eq:BPIexp-conv}
\sup_{(i,j) \in E} \bigl\| \nu_{i\to j}^{(t)} - \nu^*_{i \to j} \bigr\|
_{\mathrm{TV}}
\le A \exp(-\lambda t)  .
\end{eqnarray}
\end{theorem}

For $i_*\in V$ let $U\equiv
\mathsf{B}_{i_*}(r)$ be the ball of radius $r$ around $i_*$ in $G$,
denoting by $E_U$ its edge set, by $\partial U$ its border
(i.e., the set of its vertices at distance $r$ from $i_*$),
and for each $i\in\partial U$ let
$j(i)$ denote any one fixed neighbor of $i$ in $U$.

Our next result shows that the probability distribution
%
\begin{eqnarray}\label{eq:LocalMarg}
\quad\ \, \nu_U(\underline{x}_U)= \frac{1}{z_U}\exp\biggl\{\beta
\sum_{(i,j)\in E_U}x_ix_j+B\sum_{i\in U\setminus\partial U}x_i
\biggr\}
\prod_{i\in\partial U}\nu^*_{i\to j(i)}(x_i)  ,
\end{eqnarray}
with $\{\nu_{i\to j}^*( \cdot )\}$ the fixed point of the BP iteration
per Theorem \ref{theorem:BPConvergence},
is a good approximation for the marginal $\mu_U(  \cdot  )$
of variables $\underline{x}_U\equiv\{x_i\dvtx  i\in U\}$ under the Ising model
(\ref{eq:IsingModel}).
\begin{theorem}\label{theorem:BPCorrectness}
Assume $\beta\ge0$, $B>0$ and $G$ is a graph of finite maximal
degree~$\Delta$.
Then, there exist
finite $c=c(\beta,B,\Delta)$ and $\lambda=\lambda(\beta,B,\Delta
)>0$ such
that for any $i_* \in G$ and $U=\mathsf{B}_{i_*}(r)$, if $\mathsf{B}_{i_*}(t)$
is a
tree then
%
\begin{eqnarray}
\|\mu_U - \nu_U\|_{\mathrm{TV}} \le
\exp\{c^{r+1}-\lambda(t-r)\}  .
\end{eqnarray}
\end{theorem}
%

\subsection{Examples}

Many common random graph ensembles \cite{Rgraphs}
naturally fit our framework.

\subsubsection*{Random regular graphs} Let $G_n$ be a uniformly random graph
with degree
$k$. As $n\to\infty$, the sequence $\{G_n\}$ is obviously uniformly sparse,
and converges locally almost surely to the rooted infinite tree of
degree $k$ at every vertex.
Therefore, in this case Theorem \ref{theorem:free_energy} applies
with $P_k=1$ and $P_i=0$ for $i\neq k$.
The distributional recursion
(\ref{eqn:h_recursion}) then evolves with a deterministic sequence $h^{(t)}$
recovering the result of \cite{GerschenMon1}.

\subsubsection*{\texorpdfstring{Erd\"{o}s--Renyi graphs}{Erdos--Renyi graphs}} Let $G_n$ be a uniformly random graph
with $m=n\gamma$ edges over $n$ vertices. The sequence $\{G_n\}$ converges
locally almost surely to a Galton--Watson tree with Poisson offspring
distribution of mean $2\gamma$. This corresponds to taking
$P_k=(2\gamma)^ke^{-2\gamma}/k!$.
The same happens to classical variants of this ensemble. For instance,
one can add an edge independently for each pair $(i,j)$ with probability
$2\gamma/n$, or consider a multi-graph with $\operatorname{Poisson}(2\gamma/n)$
edges between each pair $(i,j)$.

The sequence
$\{G_n\}$ is with probability one uniformly sparse in
each of these cases. Thus,
Theorem \ref{theorem:free_energy} extends the results of
\cite{deSantisGuerra} to arbitrary nonzero temperature and magnetic field.

\subsubsection*{Arbitrary degree distribution} Let $P$ be a distribution with
finite second moment and $G_n$ a uniformly random graph with degree
distribution $P$
(more precisely, we set
the number of vertices of degree $k \ge1$ to
$\lfloor nP_k \rfloor$, adding one for $k=1$ if needed for
an even sum of degrees). Then, $\{G_n\}$ is uniformly sparse and
with probability one it
converges locally to $\mathsf{T}(P,\rho ,\infty)$.
The same happens
if $G_n$ is drawn according to the so-called configuration model
(cf. \cite{Bollobas}).

%

\section{Preliminaries}
\label{sec:Notations}

We review here the notations and a couple of classical tools
we use throughout this paper. To this end, when
proving our results it is useful to allow for vertex-dependent
magnetic fields $B_i$, that is, to replace the
basic model (\ref{eq:IsingModel}) by
%
\begin{eqnarray}\label{eq:IsingModelGen}
\mu(\underline{x}) =\frac{1}{Z(\beta,\underline{B})}  \exp\biggl\{
\beta\sum
_{(i,j)\in E}x_ix_j
+\sum_{i\in V}B_ix_i\biggr\}  .
\end{eqnarray}

Given $U\subseteq V$, we denote by $(+)_U$ [respectively, $(-)_U$]
the vector $\{x_i = +1$,  $i\in U\}$ [respectively,
$\{x_i = -1,  i\in U\}$], dropping the subscript $U$ whenever
clear from the context. Further, we use $\underline{x}_U\preceq
\underline{x}'_U$
when two real-valued vectors $\underline{x}$ and $\underline{x}'$ are
such that
$x_i\le x_i'$ for all $i\in U$ and say that
a distribution $\rho_U( \cdot )$ over ${\mathbb{R}}^U$ is
dominated by a distribution $\rho_U'( \cdot )$ over
this set (denoted $\rho_U\preceq\rho'_U$),
if the two distributions can be coupled so
that $\underline{x}_U\preceq\underline{x}'_U$ for any
pair $(\underline{x}_U,\underline{x}'_U)$ drawn from this coupling.
Finally, we use throughout the shorthand
$\langle\nu, f \rangle=\sum_x f(x) \nu(x)$
for a distribution $\nu$ and function $f$ on the same finite set, or
$\langle f \rangle$ when $\nu$ is clear from the context.

The first classical result we need is Griffiths inequality (see
\cite{Liggett}, Theorem~IV.1.21).
\begin{theorem}\label{theorem:Monotonicity}
Consider two Ising models $\mu( \cdot )$ and $\mu'( \cdot )$
on graphs $G=(V,E)$ and $G'=(V,E')$,
inverse temperatures $\beta$ and $\beta'$, and magnetic fields $\{
B_i\}$
and $\{B_i'\}$, respectively.
If $E\subseteq E'$, $\beta\le\beta'$ and $0 \le B_i \le B_i'$
for all $i\in V$, then
$0 \le\langle\mu,\prod_{i \in U}   x_i \rangle\le\langle\mu
',\prod_{i \in U} x_i \rangle$
for any $U \subseteq V$.
\end{theorem}

The second classical result we use is the GHS inequality (see \cite{GHS})
about the effect of the magnetic field $\underline{B}$ on the local
magnetizations
at various vertices.
\begin{theorem}[(Griffiths, Hurst, Sherman)]\label{theorem:GHS}
Let $\beta\ge0$ and for $\underline{B}= \{B_i\dvtx  i\in V\}$, denote by
$m_j(\underline{B})\equiv\mu(\{\underline{x}\dvtx x_j=+1\})-\mu(\{
\underline{x}\dvtx x_j=-1\})$
the local magnetization at vertex $j$ in
the Ising model (\ref{eq:IsingModelGen}).
If $B_i\ge0$ for all $i\in V$, then for any three
vertices $j,k,l\in V$ (not necessarily distinct),
%
\begin{eqnarray}
\frac{\partial^2 m_j(\underline{B})}{\partial B_k\,\partial B_l}\le
0  .
\end{eqnarray}
\end{theorem}

Finally, we need the following elementary inequality:
\begin{lemma}\label{lemma:SimpleIneq}
For any function $f\dvtx \mathcal{X}\mapsto[0,f_{\max}]$
and distributions $\nu$, $\nu'$ on the finite
set $\mathcal{X}$ such that $\nu(f>0)>0$ and $\nu'(f>0)>0$,
%
\begin{eqnarray}\label{eq:SimpleIneq}
\sum_x\biggl|\frac{\nu(x)f(x)}{\langle\nu,f\rangle}-\frac{\nu
'(x)f(x)}{\langle\nu',f \rangle}
\biggr|\le\frac{3f_{\max}}{\max(\langle\nu,f \rangle,\langle
\nu',f \rangle)} \|\nu-\nu'\|_{\mathrm{TV}}.
\end{eqnarray}
In particular, if $0<f_{\min}\le f(x)$, then the
right-hand side is bounded by $(3f_{\max}/f_{\min}) \|\nu-\nu'\|_{\mathrm{TV}}$.
\end{lemma}

\begin{pf}
Assuming without loss of generality that $\langle\nu',f \rangle\ge
\langle\nu,f\rangle>0$, the left-hand side of (\ref{eq:SimpleIneq}) can be bounded as
\begin{eqnarray*}
&&\frac{1}{\langle\nu,f \rangle\langle\nu',f \rangle}
\sum_x|\nu(x)f(x) \langle\nu',f \rangle-\nu'(x)f(x) \langle
\nu,f \rangle|
\\
&&\qquad \le\frac{1}{\langle\nu',f\rangle}|\langle\nu,f \rangle-
\langle\nu',f \rangle|+
\frac{1}{\langle\nu',f\rangle}\sum_x|\nu(x)f(x)-\nu'(x)f(x)|
\\
&&\qquad\le\frac{f_{\max}}{\langle\nu',f\rangle}\|\nu-\nu'\|_{\mathrm{TV}}+
\frac{2f_{\max}}{\langle\nu',f \rangle} \|\nu-\nu'\|_{\mathrm
{TV}}  .
\end{eqnarray*}
This implies the lemma.
\end{pf}
%

%
%
\section{Ising models on trees}
\label{sec:TreeProofs}

We prove in this section certain facts about Ising models on trees
which are of independent interest and as a byproduct we
deduce Lemma \ref{lemma:Recursive} and the theorems
of Section \ref{sec:AlgoImpl}.
In doing so, recall that
for each $\ell\ge1$ the Ising models
on $\mathsf{T}(\ell)$ with free and plus boundary conditions are
\begin{eqnarray}
\mu^{\ell,0}(\underline{x}) &\equiv&\frac{1}{Z^{\ell,0}}
\exp\biggl\{\beta\sum_{(ij)\in\mathsf{T}(\ell)}x_ix_j+\sum_{i\in
\mathsf{T}(\ell)}
B_ix_i\biggr\}  ,
\\
\mu^{\ell,+}(\underline{x}) &\equiv&\frac{1}{Z^{\ell,+}}
\exp\biggl\{\beta\sum_{(ij)\in\mathsf{T}(\ell)}x_ix_j+\sum_{i\in
\mathsf{T}(\ell)}
B_ix_i\biggr\}\nonumber
\\[-8pt]\\[-8pt]
&&{}\times \mathbb{I}\bigl(\underline{x}_{\partial\mathsf{T}(\ell
)} = (+)_{\partial
\mathsf{T}(\ell)}\bigr)  .\nonumber
\end{eqnarray}
Equivalently $\mu^{\ell,0}$ is the Ising model (\ref{eq:IsingModelGen})
on $\mathsf{T}(\ell)$
with magnetic fields\vspace{1pt} $\{B_i\}$ and $\mu^{\ell,+}$ is the modified Ising
model corresponding to the limit $B_i \uparrow+\infty$
for all $i\in\partial\mathsf{T}(\ell)$. To simplify our notation we
denote such limits hereafter simply by setting $B_i=+\infty$
and use $\mu^{\ell}$ for statements that apply to both free
and plus boundary conditions.

We start with the following simple but useful observation.
\begin{lemma}\label{lemma-observation}
For a subtree $U$ of a finite tree $T$ let
$\partial_* U$ denote the subset of
vertices of $U$ connected by an edge
to $W \equiv T \setminus U$ and
for each $u \in\partial_* U$ let
$\langle x_u\rangle_W$ denote the root magnetization of the Ising model
on the maximal subtree $T_u$ of $W \cup\{u\}$
rooted at $u$. The marginal on $U$ of the Ising measure on $T$,
denoted $\mu^T_U$ is then an Ising measure on $U$ with
magnetic field $B_u' = \operatorname{atanh}(\langle x_u\rangle_W) \ge B_u$ for
$u \in\partial_* U$ and $B_u' = B_u$ for $u \notin\partial_* U$.
\end{lemma}
\begin{pf} Since $U$ is a subtree
of the tree $T$, the subtrees
$T_u$ for $u \in\partial_* U$ are disjoint.
Therefore, with $\hat{\mu}_u(\underline{x})$ denoting
the Ising model distribution for $T_u$
we have that
%
\begin{eqnarray}\label{eq:boltz}
\mu_{U}^T(\underline{x}_U) = \frac{1}{\hat{Z}}
f(\underline{x}_U)
\prod_{u\in\partial_* U}\hat{\mu}_u(x_u)
\end{eqnarray}
for the Boltzmann weight
\[
f(\underline{x}_U)=\exp\biggl\{\beta\sum_{(uv)\in U}x_u x_v
+\sum_{u \in U \setminus\partial_* U} B_u x_u \biggr\}  .
\]
Further, $x_u \in\{+1,-1\}$ so for each $u \in\partial_* U$
and some constants $c_u$,
\[
\hat{\mu}_u (x_u) = \tfrac{1}{2}(1+ x_u \langle x_u\rangle_W) = c_u
\exp(\operatorname{atanh}(\langle x_u\rangle_W) x_u)  .
\]
Embedding the normalization constants $c_u$ within $\hat{Z}$
we thus conclude that $\mu^T_U$ is an Ising measure on $U$
with the stated magnetic field $B_u'$. Finally, comparing
the root magnetization for $T_u$ with that for
$\{u\}$ we have by Griffiths
inequality that $\langle x_u\rangle_W \ge\tanh(B_u)$, as claimed.
\end{pf}

\begin{theorem}\label{theorem:TreeDecay}
Suppose $\mathsf{T}$ is a conditionally independent infinite tree
of average offspring numbers bounded by $\Delta$,
as in Definition \ref{def-cit}.
For $0 < B_{\min} \le B_{\max}$, $\beta_{\max}$ and $\Delta$ finite,
there exist $M=M(\beta_{\max},B_{\min},\Delta)$ and
$C=C(\beta_{\max},B_{\max})$ finite such that if
$B_i \le B_{\max}$ for all $i \in\mathsf{T}(r-1)$ and
$B_i\ge B_{\min}$ for all $i \in\mathsf{T}(\ell)$, $\ell> r$, then
%
\begin{eqnarray}\label{eq:bd-influence}
\mathbb{E}  \|\mu_U^{\ell,+}-\mu_U^{\ell,0}\|_{\mathrm{TV}}\le
\delta(\ell-r)   \mathbb{E}\bigl\{C^{|\mathsf{T}(r)|}\bigr\}
\end{eqnarray}
for $\delta(t)=M/t$, all $U \subseteq\mathsf{T}(r)$ and $\beta\le
\beta
_{\max}$.
\end{theorem}

%
\begin{pf} Fixing $\ell> r$
it suffices to consider $U=\mathsf{T}(r)$ [for which
the \mbox{left-hand} side of (\ref{eq:bd-influence}) is maximal].
For this $U$ and $T=\mathsf{T}(\ell)$ we have that
$\partial_* U=\partial\mathsf{T}(r)$ and
$U \setminus\partial_* U = \mathsf{T}(r-1)$, where in
this case the Boltzmann weight $f(\cdot)$ in (\ref{eq:boltz})
is bounded above by $f_{\max}=c^{|\mathsf{T}(r)|}$ and
below by $f_{\min}=1/f_{\max}$
for $c=\exp(\beta_{\max}+B_{\max})$. Further, the
plus and free boundary conditions then differ in
(\ref{eq:boltz}) by having the corresponding
boundary conditions at generation $\ell-r$ of
each subtree $T_u$, which we
distinguish by using
$\hat{\mu}^{+/0}_u(x_u)$ instead of
$\hat{\mu}_u(x_u)$. Since
the total variation distance between two product measures
is at most the sum of the distance between their marginals,
upon applying Lemma \ref{lemma:SimpleIneq} we deduce
from (\ref{eq:boltz}) that
\[
\bigl\| \mu_{\mathsf{T}(r)}^{\ell,+} - \mu_{\mathsf{T}(r)}^{\ell,0}
\bigr\|_{\mathrm{TV}}
\le\frac{3}{2} c^{2 |\mathsf{T}(r)|}
\sum_{i \in\partial\mathsf{T}(r)} | \hat{\mu}^{+}_i(x_i=1)
- \hat{\mu}^{0}_i(x_i=1) |  .
\]
By our assumptions, conditional on $U=\mathsf{T}(r)$,
the subtrees $T_i$ of $T=\mathsf{T}(\ell)$
denoted hereafter also by $\mathsf{T}_i$ are
for $i \in\partial\mathsf{T}(r)$ independent of each other.
Further, $2 \hat{\mu}^{+/0}_i(x_i=1)-1$
is precisely the magnetization of their root vertex
under plus/free boundary conditions at generation
$\ell-r$. Thus, taking $C=e c^2$
(and using the inequality
$y \le e^y$), it suffices to show that the magnetizations
$m^{\ell,+/0}(\underline{B}) = \langle\mu^{\ell,+/0},x_{\o}
\rangle$
at the root of any such conditionally independent infinite tree
$\mathsf{T}$ satisfy
$\mathbb{E}\{m^{\ell,+}(\underline{B})-m^{\ell,0}(\underline{B})\}
\le\frac{M}{\ell}$,
for some $M=M(\beta_{\max},B_{\min},\Delta)$ finite,
all $\beta\le\beta_{\max}$ and
$\ell\ge1$, where we have removed the absolute value since
$m^{\ell,+}(\underline{B})\ge m^{\ell,0}(\underline{B})$ by
Griffiths inequality.
For greater convenience of the reader, this fact is
proved in the next lemma.
\end{pf}

\begin{lemma}\label{lemma:mg-root}
Suppose $\mathsf{T}$ is a conditionally independent infinite tree
of average offspring numbers bounded by $\Delta$.
For $0 < B_{\min} \le B_{\max}$, $\beta_{\max}$ and $\Delta$ finite,
there exist $M=M(\beta_{\max},B_{\min},\Delta)$ such that
%
\begin{eqnarray}
\label{eq:mg-root}
\mathbb{E}\{m^{\ell,+}(\underline{B})-m^{\ell,0}(\underline{B})\}
\le\frac{M}{\ell}   ,
\end{eqnarray}
where $m^{\ell,+/0}(\underline{B}) = \langle\mu^{\ell,+/0},x_{\o}
\rangle$
are the root magnetizations under $+$ and free boundary
condition on $\mathsf{T}$.
\end{lemma}

\begin{pf}
Note that (\ref{eq:mg-root}) trivially holds for $\beta=0$
[in which case $\mu^{\ell,+}(x_{\o})=\mu^{\ell,0}(x_{\o})$].
Assuming hereafter that $\beta>0$ we
proceed to prove (\ref{eq:mg-root})
when each vertex of $\mathsf{T}(\ell-1)$ has a
nonzero offspring number.
To this end,
for $\underline{H}=\{H_i\in{\mathbb{R}}\dvtx i\in\partial\mathsf
{T}(k)\}$ let
\begin{eqnarray*}\label{eq:muBH}
\mu^{k,\underline{H}}(\underline{x}) &\equiv&\frac{1}{Z^{k,0}}
\exp\biggl\{\beta\sum_{(ij)\in\mathsf{T}(k)}x_ix_j+\sum_{i\in
\mathsf{T}(k)}
B_ix_i+ \sum_{i\in\partial\mathsf{T}(k)}H_ix_i\biggr\}
\end{eqnarray*}
and denote by $m^{k}(\underline{B},\underline{H})$ the corresponding
root magnetization.
Writing $H$ instead of $\underline{H}$ for constant magnetic field
on the leave nodes, that is, when $H_i=H$
for each $i\in\partial\mathsf{T}(k)$, we note that
$m^{k,+}(\underline{B}) = m^{k}(\underline{B},\infty)$
and $m^{k,0}(\underline{B}) = m^{k}(\underline{B},0)$.
Further, applying Lemma \ref{lemma-observation} for
the subtree $\mathsf{T}(k-1)$ of $\mathsf{T}(k)$
we represent $m^{k}(\underline{B},\infty)$ as the root magnetization
$m^{k-1}(\underline{B}',0)$ on $\mathsf{T}(k-1)$ where
$B'_i = B_i +\beta\Delta_i$ for $i\in\partial\mathsf{T}(k-1)$ and
$B'_i = B_i$ for all other $i$. Consequently,
%
\begin{eqnarray}
\label{eq:bdpp}
m^{k}(\underline{B},\infty) = m^{k-1}(\underline{B},\{\beta\Delta
_i\})  .
\end{eqnarray}
Recall that if
$\frac{\partial^2 g}{\partial^2 z_i}\le0$ for $i=1,\ldots,s$,
then applying Jensen's inequality one variable at a time we
have that
$\mathbb{E}  g(Z_1,\dots,Z_s)\le g(\mathbb{E} Z_1,\dots,\mathbb
{E}  Z_s)$
for any independent random variables $Z_1,\ldots, Z_s$.
By the GHS inequality, this is the case for
$\underline{H}\mapsto m^{k-1}(\underline{B},\underline{H})$, hence
with $\mathbb{E}_k$ denoting
the conditional on $\mathsf{T}(k)$
expectation over the independent offspring numbers
$\Delta_i$ for $i \in\partial\mathsf{T}(k)$, we deduce that
%
\begin{eqnarray}\label{eq:bdp}
\mathbb{E}_{k-1}  m^{k}(\underline{B},\infty)
\le m^{k-1}(\underline{B},\{\beta\mathbb{E}_{k-1} \Delta_i\})
\le m^{k-1}(\underline{B},\beta\Delta)  ,
\end{eqnarray}
where the last inequality is a consequence of
Griffiths inequality and our assumption that
$\mathbb{E}_t \Delta_i \le\Delta$ for any $i \in\partial\mathsf{T}(t)$
and all $t \ge0$.
Since each $i \in\partial\mathsf{T}(k-1)$ has at least
one offspring whose magnetic field is at least $B_{\min}$,
it follows by Griffiths inequality that
$m^{k,0}(\underline{B})$ is bounded below by the
magnetization at the root of the subtree $T$ of
$\mathsf{T}(k)$ where $\Delta_i=1$ for all
$i \in\partial\mathsf{T}(k-1)$
and $B_i=B_{\min}$ for all $i \in\partial\mathsf{T}(k)$.
Applying Lemma \ref{lemma-observation}
for $T$ and $U=\mathsf{T}(k-1)$, the root magnetization for
the Ising distribution on $T$ turns out to be
precisely $m^{k-1}(\underline{B},\xi)$ for $\xi=\xi(\beta,B_{\min})>0$
of (\ref{eq:xidef}).
Thus, one more application of Griffiths inequality yields that
%
\begin{eqnarray}\label{eq:bd0}
m^{k}(\underline{B},0)\ge m^{k-1}(\underline{B},\xi)\ge
m^{k-1}(\underline{B},0)  .
\end{eqnarray}
Next note that $\xi(\beta,B) \le\beta\le\beta\Delta$ and by
GHS inequality $H \mapsto m^{k-1}(\underline{B},H)$ is concave. Hence,
%
\begin{eqnarray}\label{eq:bd1}
m^{k-1}(\underline{B},\beta\Delta) - m^{k-1}(\underline{B},0)
\le M [ m^{k-1}(\underline{B},\xi) - m^{k-1}(\underline{B},0) ]
\end{eqnarray}
for the finite constant
%
\[
M \equiv\sup_{0 < \beta\le\beta_{\max}}
  \frac{\beta\Delta}{\xi(\beta,B_{\min})}
\]
and all $\beta\le\beta_{\max}$.
Combining (\ref{eq:bdp}), (\ref{eq:bd0}) and
(\ref{eq:bd1}) we obtain that
\begin{eqnarray*}
\mathbb{E}_{k-1}\{m^{k,+}(\underline{B})-m^{k,0}(\underline{B})\}
&\le& m^{k-1}(\underline{B},\Delta\beta)- m^{k-1}(\underline{B},0)
\\
&\le& M [m^{k-1}(\underline{B},\xi)- m^{k-1}(\underline{B},0)]
\\
&\le&
M [m^{k}(\underline{B},0)- m^{k-1}(\underline{B},0)]   .
\end{eqnarray*}
We have seen in (\ref{eq:bd0}) that
$k \mapsto m^{k,0}(\underline{B})$ is nondecreasing whereas
from (\ref{eq:bdpp}) and Griffiths inequality we have
that $k \mapsto m^{k,+}(\underline{B})$ is nonincreasing.
With magnetization bounded above by one, we thus get
upon summing the preceding inequalities for $k=1,\ldots,\ell$ that
\begin{eqnarray*}
\ell\mathbb{E}_{\ell-1} [ m^{\ell,+}(\underline{B})-m^{\ell
,0}(\underline{B}) ]\le
\sum_{k=1}^{\ell} \mathbb{E}_{k-1} [ m^{k,+}(\underline
{B})-m^{k,0}(\underline{B}) ]\le M   ,
\end{eqnarray*}
from which we deduce (\ref{eq:mg-root}).

Considering now the general case where the
infinite tree $\mathsf{T}$ has vertices (other than the root)
of degree one, let $\mathsf{T}^*(\ell)$ denote the ``backbone''
of $\mathsf{T}(\ell)$, that is, the subtree
induced by vertices along self-avoiding paths between $\o$ and
$\partial\mathsf{T}(\ell)$. Taking $U=\mathsf{T}^*(\ell)$ as the subtree
of $T=\mathsf{T}(\ell)$ in Lemma \ref{lemma-observation},
note that for each $u \in\partial_* U$ the
subtree $T_u$ contains no vertex from $\partial\mathsf{T}(\ell)$.
Consequently, the marginal measures $\mu^{\ell,+/0}_U$ are Ising
measures on $U$ with the same magnetic fields
$B_i' \ge B_i \ge B_{\min}$ outside $\partial\mathsf{T}(\ell)$.
Thus, with $m^{\ell,+/0}_*(\underline{B})$
denoting the corresponding magnetizations at the root
for $\mathsf{T}^*(\ell)$, we deduce that
$m^{\ell,+/0}(\underline{B})=m^{\ell,+/0}_*(\underline{B}')$ where
$B'_i \ge B_i \ge B_{\min}$ for all $i$. By definition every
vertex of $\mathsf{T}^*(\ell-1)$ has a nonzero
offspring number and with $B_i' \ge B_{\min}$, the required bound
\[
\mathbb{E}\{m^{\ell,+}(\underline{B})-m^{\ell,0}(\underline{B})\} =
\mathbb{E}\{m^{\ell,+}_*(\underline{B}')-m^{\ell,0}_*(\underline
{B}')\} \le\frac{M}{\ell}
\]
follows by the preceding argument, since
$\mathsf{T}^*(\ell)$ is a conditionally independent tree
whose offspring numbers $\Delta_i^* \ge1$ do not exceed those of
$\mathsf{T}(\ell)$. Indeed, for $k=0,1,\ldots, \ell-1$,
given $\mathsf{T}^*(k)$
the offspring numbers at $i \in\partial\mathsf{T}^*(k)$
are independent of each other [with probability of
$\{\Delta^*_i=s\}$ proportional to the sum over $t \ge0$
of the product of the probability
of $\{\Delta_i=s+t\}$ and that of precisely $s$ out of the $s+t$
offspring of $i$ in $\mathsf{T}(\ell)$ having a line of descendants
that survives additional $\ell-k-1$ generations, for $s \ge1$].
\end{pf}

Simon's inequality (see \cite{Simon}, Theorem 2.1)
allows one to bound the (centered) two point
correlation functions in ferromagnetic
Ising models with zero magnetic field.
We provide next its generalization to arbitrary magnetic
field, in the case of Ising models on trees.
\begin{lemma}\label{eq:LemmaPatch}
If edge $(i,j)$ is on the
unique path from $\o$ to $k\in\mathsf{T}(\ell)$, with
$j$ a descendant of $i \in\partial\mathsf{T}(t)$, $t \ge0$, then
%
\begin{eqnarray}\label{ineq-simon}
\langle x_{\o};x_k\rangle^{(\ell)}_{\o}\le
\cosh^2(2\beta+B_i)
  \langle x_{\o};x_{i}\rangle^{(t)}_{\o}
\langle x_{j};x_{k}\rangle^{(\ell)}_{j}  ,
\end{eqnarray}
where $\langle \cdot \rangle^{(r)}_{i}$ denotes the expectation
with respect to
the Ising distribution $\hat{\mu}_i(\cdot)$
on the subtree $\mathsf{T}_i$ of $i$ and all its descendants in
$\mathsf{T}(r)$ and $\langle x;y\rangle\equiv\langle xy\rangle
-\langle x\rangle\langle y\rangle$ denotes the
centered two point correlation function.
\end{lemma}
\begin{pf}
It is not hard to
check that if $x,y,z$ are $\{+1,-1\}$-valued random variables
with $x$ and $z$ conditionally independent given $y$, then
%
\begin{eqnarray}
\langle x;z\rangle= \frac{\langle x;y\rangle\langle y;z\rangle
}{1-\langle y\rangle^2}   .
\label{eq:Identity}
\end{eqnarray}
In particular, under $\mu^{\ell,0}$ the random variables
$x_{\o}$ and $x_k$ are conditionally independent given
$y=x_i$ with
\[
\biggl| \log\biggl( \frac{\mu^{\ell,0}(x_i=+1)}{\mu^{\ell
,0}(x_i=-1)} \biggr)
\biggr|
\le2 (|\partial i| \beta+B_i)  .
\]
Hence, if $j$ is the unique descendant of $i$
then $|\langle x_i\rangle_{\o}^{(\ell)}| \le\tanh(2\beta+B_i)$
and we get from (\ref{eq:Identity}) that
\[
\langle x_{\o};x_k\rangle^{(\ell)}_{\o} \le c
\langle x_{\o};x_i\rangle^{(\ell)}_{\o}\langle x_i;x_k\rangle
^{(\ell)}_{\o}
\]
for $c=\cosh^2(2\beta+B_i)$.
Next note that $\langle x;y\rangle\le1-\langle y\rangle^2$ for any two
$\{+1,-1\}$-valued random variables, and since
$x_i$ and $x_k$ are conditionally independent given
$y=x_j$ it follows from (\ref{eq:Identity}) that
$\langle x_i;x_k\rangle^{(\ell)}_{\o}\le\langle x_j;x_k\rangle
^{(\ell)}_{\o}$.
Further, if $\langle\cdot\rangle$ is the expectation with respect to
an Ising measure for some (finite) graph $G$ then
for any $u,v \in G$
%
\begin{eqnarray}\label{eq:ident-corr}
\frac{\partial\langle x_v\rangle}{\partial B_u} =
\langle x_v x_u\rangle-\langle x_v\rangle \langle x_u\rangle =
\langle x_v; x_u \rangle  .
\end{eqnarray}
From Lemma \ref{lemma-observation} we
know that computing the marginal of the Ising distribution
for $T=\mathsf{T}(\ell)$ on a smaller subtree $U=\mathsf{T}_j$ of interest
has the effect of increasing its magnetic field.
Thus, combining the identity (\ref{eq:ident-corr}) with
GHS inequality, we see that
reducing this field (i.e., restricting to $U$ the original
Ising distribution), increases the
centered two point correlation function. That is,
$\langle x_j;x_k\rangle^{(\ell)}_{\o}\le\langle x_j;x_k\rangle
^{(\ell)}_{j}$.
Similarly, considering Lemma \ref{lemma-observation}
for $U=\mathsf{T}(t)$ we also\vspace{1pt} have that
$\langle x_{\o};x_i\rangle^{(\ell)}_{\o}\le
\langle x_{\o};x_i\rangle^{(t)}_{\o}$ which completes our thesis
in case $j$ is the unique descendant of $i$.

Turning to the general case, we compare the thesis of the lemma for
$\mathsf{T}(\ell)$ and the subtree $U=\mathsf{T}'(\ell)$ obtained upon
deleting the subtrees rooted at descendants of $i$
(and the corresponding edges to $i$) except for $\mathsf{T}_j$.
While $\langle x_{\o};x_{i}\rangle^{(t)}_{\o}$ and
$\langle x_{j};x_{k}\rangle^{(\ell)}_{j}$ are unchanged by this modification
of the underlying tree (as the relevant subgraphs are not modified),
we have from Lemma \ref{lemma-observation} that
$\mu_U^{\ell,0}(\cdot)$ is an Ising measure on $U$
identical to the original but for an increase in the
magnetic field at $i$. In view of
(\ref{eq:ident-corr}) and the GHS inequality, we thus deduce
that the value of $\langle x_{\o};x_k\rangle^{(\ell)}_{\o}$
is smaller for the Ising model on $\mathsf{T}(\ell)$ than
for the one on $\mathsf{T}'(\ell)$ and since in $\mathsf{T}'(\ell)$ the
vertex $j$ is the unique descendant of $i$, we are done.
\end{pf}

Equipped with the preceding lemma we
next establish the exponential decay of correlations
and of the effect of boundary conditions in Theorem
\ref{theorem:TreeDecay}.
\begin{coro}\label{coro:Susc}
There exist $A$ finite and $\lambda$ positive, depending
only on $\beta_{\max}$, $B_{\min}$, $B_{\max}$
and $\Delta$ such that
%
\begin{eqnarray}\label{eq:ExponentialDecay}
\mathbb{E}\biggl\{\sum_{i\in\partial\mathsf{T}(r)}  \langle x_{\o
};x_i\rangle_{\o
}^{(\ell)}
\biggr\}\le A   e^{-\lambda r}
\end{eqnarray}
for any $r \le\ell$ and if $B_i \le B_{\max}$ for all $i \in\mathsf{T}
(\ell-1)$
then Theorem \ref{theorem:TreeDecay} holds for $\delta(t) = A \exp
(-\lambda t)$.
\end{coro}

\begin{Remark*}
Taking $B_i \uparrow+\infty$ for $i \in\partial\mathsf{T}(\ell)$,
note that (\ref{eq:ExponentialDecay}) applies
when $\langle \cdot \rangle^{(\ell)}$ is with respect to
$\mu^{\ell,+}( \cdot )$.
\end{Remark*}

\begin{pf*}{Proof of Corollary \ref{coro:Susc}} Starting with the proof of (\ref{eq:ExponentialDecay})
take $\ell= r$ for which the left-hand side is maximal
(as we have seen while proving Lemma \ref{eq:LemmaPatch}).
Then, denoting by $\langle  \cdot  \rangle_{H_r}$ the
expectation under the Ising measure on $\mathsf{T}(r)$ with
a magnetic field $H_r$ added to $\underline{B}$ at all vertices
$i \in\partial\mathsf{T}(r)$, it follows from
(\ref{eq:ident-corr}) that
\begin{eqnarray*}
\sum_{i\in\partial\mathsf{T}(r)}   \langle x_{\o}; x_i\rangle
_{\o}^{(r)}
= \sum_{i\in\partial\mathsf{T}(r)}
\frac{\partial\langle x_{\o}\rangle}{\partial B_i} =
\frac{\partial\langle x_{\o}\rangle_{H_r}}{\partial H_r}
\bigg|_{H_r=0}  .
\end{eqnarray*}
By GHS inequality the latter derivative is nonincreasing
in $H_r$, whence
\begin{eqnarray*}
\sum_{i\in\partial\mathsf{T}(r)}   \langle x_{\o}; x_i\rangle
_{\o}^{(r)}
\le \frac{2}{B_{\min}}[\langle x_{\o}\rangle_{H_r=0}-\langle
x_{\o}\rangle_{
H_r=-B_{\min}/2}]  .
\end{eqnarray*}
Let $B_i'=B_i-B_{\min}/2$ if
$i \in\partial\mathsf{T}(r)$ and $B_i'=B_i$ otherwise, so
$\langle x_{\o}\rangle_{H_r=-B_{\rm min}/2} = m^{r,0}(\underline{B}')$.
Further, from Griffiths inequality also
$\langle x_{\o}\rangle_{H_r=0}\le\langle x_{\o}\rangle_{H_r=\infty}
= m^{r,+}(\underline{B}')$ and it follows that
%
\begin{eqnarray}\label{eq:bd-corr}
\Gamma_{r} \equiv
\mathbb{E}\biggl\{\sum_{i\in\partial\mathsf{T}(r)}   \langle x_{\o
};x_i\rangle_{\o}^{(r)}
  \biggr\}
\le \frac{2}{B_{\min}}\mathbb{E}\{m^{r,+}(\underline
{B}')-m^{r,0}(\underline{B}')\}  .
\end{eqnarray}
In particular, setting $c=\cosh^2(2\beta_{\max}+B_{\max})$, in view
of Lemma \ref{lemma:mg-root} we find that
$\Gamma_{d - 1} \le1/(e c \Delta)$ for
$d = 1 + \lceil2 e c \Delta
M(\beta_{\max},B_{\min}/2,\Delta)/B_{\min} \rceil$.
Further, since $\mathsf{T}$ is conditionally independent,
the same proof shows that if $t + d = r' \le r$ and
$\mathsf{T}_j$ is the subtree of $\mathsf{T}(r)$ of depth $d - 1$
rooted at $j \in\partial\mathsf{T}(t+1)$ then
\[
\mathbb{E}_{t+1} \biggl\{\sum_{k\in\partial\mathsf{T}_j}   \langle
x_j;x_k\rangle
_{j}^{(r')}
\biggr\} \le\frac{1}{e c \Delta}  .
\]
Considering  inequality (\ref{ineq-simon}) of Lemma \ref{eq:LemmaPatch}
for $t = r -d \equiv r_1$
and all $k \in\partial\mathsf{T}(r)$ we find that
\begin{eqnarray*}
\Gamma_{r} &\le& c
\mathbb{E}\biggl\{ \mathop{\sum_{
i\in\partial\mathsf{T}(t)}}_{
j\in\partial\mathsf{T}(t+1)\cap\partial i}
\langle x_{\o};x_{i}\rangle^{(t)}_{\o} \mathbb{E}_{t+1} \biggl[
\sum_{k\in\partial\mathsf{T}_j}\langle x_{j};x_{k}\rangle
^{(r)}_{j} \biggr] \biggr\}
\\
&\le& \frac{1}{e \Delta}
\mathbb{E}\biggl\{ \sum_{i\in\partial\mathsf{T}(t)}
\Delta_i \langle x_{\o};x_{i}\rangle^{(t)}_{\o}
\biggr\} \le e^{-1} \Gamma_{r_1}  .
\end{eqnarray*}
Iterating the preceding bound at $r_s = r - s d$, for
$s=1,\ldots, \lfloor r/d \rfloor$ and noting that
by (\ref{eq:bd-corr}) we have the bound $\Gamma_{r'} \le2/B_{\min}$
at the last step, we get the uniform in $\beta\le\beta_{\max}$
exponential decay of (\ref{eq:ExponentialDecay}).

Next, recall that the rate $\delta(t)$
in Theorem \ref{theorem:TreeDecay} is merely
the rate in the bound (\ref{eq:mg-root}).
For $k \equiv|\partial\mathsf{T}(\ell)|$ we
choose uniformly and independently of everything else a
one to one mapping
$i\dvtx \{1,\ldots,k\} \mapsto\partial\mathsf{T}(\ell)$,
and let $\underline{B}^{(s)}$ for $s \ge1$ denote
the magnetic field configuration obtained when taking
$B_{i(j)} \uparrow+\infty$ for all $j \le s$
(with $\underline{B}^{(0)} = \underline{B}$). Since
\[
m^{\ell,+}(\underline{B})-m^{\ell,0}(\underline{B}) = \sum_{s=0}^{k-1}
\bigl[m^{\ell,0}\bigl(\underline{B}^{(s+1)}\bigr)-m^{\ell,0}\bigl(\underline{B}^{(s)}\bigr)
\bigr]  ,
\]
we get the rate $\delta(t) = A \exp(-\lambda t)$
from (\ref{eq:ExponentialDecay}) as soon as we show that
for $i=i(s+1)$ and $s=0,\ldots,k-1$,
%
\begin{eqnarray}\label{eq:MagnVsSusc}
m^{\ell,0}\bigl(\underline{B}^{(s+1)}\bigr)-m^{\ell,0}\bigl(\underline{B}^{(s)}\bigr)
\le
\langle x_{\o};x_i\rangle_{\o}^{(\ell)}  .
\end{eqnarray}
To this end, let $\langle \cdot \rangle_s$ denote
the expectation under $\mu^{\ell,0}$ with
magnetic field $\underline{B}^{(s)}$
so $m^{\ell,0}(\underline{B}^{(s)})=\langle x_{\o}\rangle_s$.
Further, fixing $i=i(s+1)$
\begin{eqnarray*}
m^{\ell,0}\bigl(\underline{B}^{(s+1)}\bigr)=
\frac{\langle x_{\o}\mathbb{I}(x_i=1)\rangle_{s}}{\langle\mathbb
{I}(x_{i}=1)\rangle_{s}}
=\frac{\langle x_{\o}x_i\rangle_s+\langle x_{\o}\rangle
_s}{1+\langle x_i\rangle_s}
\end{eqnarray*}
[since $\mathbb{I}(x_i=1)=(1+x_i)/2$]. Since
$\langle x_i\rangle_s \ge0$ by Griffiths inequality, it follows
that
\[
m^{\ell,0}\bigl(\underline{B}^{(s+1)}\bigr)-m^{\ell,0}\bigl(\underline{B}^{(s)}\bigr)
\le
\langle x_{\o}x_i\rangle_s-\langle x_{\o}\rangle_s\langle
x_i\rangle_s =
\frac{\partial m_{\o}(\underline{B}^{(s)})}{\partial B_i}
  ,
\]
which by GHS inequality is maximal at $s=0$, yielding
(\ref{eq:MagnVsSusc}) and completing the proof.
\end{pf*}

As promised, Lemma \ref{lemma:Recursive} follows from the preceding results.

\begin{pf*}{Proof of Lemma \ref{lemma:Recursive}}
Consider the Galton--Watson tree $\mathsf{T}(\rho ,\infty)$
of Section~\ref{sec:Local} and the corresponding
Ising models $\mu^{t,+/0}(\underline{x})$
of constant magnetic field $B_i = B>0$ on
the subtrees $\mathsf{T}(\rho ,t)$.
It is easy to check that the random variables
$h^{(t)} = \operatorname{atanh}(m^{t,0}(B))$
satisfy the distributional recursion (\ref{eqn:h_recursion})
starting at $h^{(0)}=0$.
By Griffiths inequality $m^{t,0}(B)$, hence $h^{(t)}$,
is nondecreasing in $t$, and so
converges almost surely as $t \to\infty$ to a limiting
random variable $h^*$. Further, the bounds
$0= h^{(0)} \le h^{(t)} \le B + \Delta_{\o}$ hold for all $t$
and hence also for $h^*$. We thus deduce that
the distributions $Q_t$ of $h^{(t)}$ as determined
by (\ref{eqn:h_recursion}) are stochastically monotone (in~$t$)
and converge weakly to some law $Q^*$ of $h^*$ that
is supported on $[0,\infty)$.

Next, recall that for any fixed $k$ and
$F(\cdot)$ continuous and bounded on ${\mathbb{R}}^k$, the functional
$\Psi_F(Q) = \int F(h_1,\ldots,h_k) \,\mathrm{d}Q(h_1) \,\cdots\,\mathrm
{d}Q(h_k)$
is continuous with respect to weak convergence
of probability measures on $[0,\infty)$ (e.g., see
\cite{DZ}, Lemma~7.3.12). Fixing $g\dvtx {\mathbb{R}}\mapsto[-C,C]$
continuous, clearly
\[
g_j (h_1,\ldots,h_j)= g\Biggl( B + \sum_{i=1}^{j-1} \xi(\beta,h_i) \Biggr)
\]
are continuous and bounded. Further,
it follows from (\ref{eqn:h_recursion}) that for all $t$
\[
\Biggl|\int g\, \mathrm{d}Q_{t+1} -
\sum_{j=1}^k \mathbb{P}(K=j)
\Psi_{g_j}(Q_t)\Biggr| \le C \mathbb{P}(K > k)  .
\]
Taking $t \to\infty$ followed by
$k \to\infty$, we deduce by the preceding arguments
[and the uniform boundedness
$|\Psi_{g_j}(Q^*)| \le C$ for all $j$], that
\[
\int g \,\mathrm{d}Q^{*} = \sum_{j=1}^\infty\mathbb{P}(K=j) \Psi
_{g_j}(Q^*)  .
\]
As this applies for every bounded continuous function
$g(\cdot)$, we conclude that $h^*$ and its law $Q^*$
are a fixed point of
the distributional recursion (\ref{eqn:h_recursion}).

Next note that the random variables
$h^{(t)}_+ = \operatorname{atanh}[m^{t,+}(B)]$ form a non-increasing
sequence that satisfies the same
distributional recursion, but with the initial condition
$h^{(0)}_+ = +\infty$. Consequently, by the same
arguments we have used before, the
laws $Q_{t,+}$ of $h^{(t)}_+$ converge weakly to
some fixed point $Q^*_+$ of (\ref{eqn:h_recursion})
that is also supported on $[0,\infty)$.
Further, $Q_t\preceq Q^{**}\preceq Q_{t,+}$
for $t=0$ and any (other) possible law $Q^{**}$
of a fixed point $h^{**}$ of (\ref{eqn:h_recursion})
that is supported on $[0,\infty)$. Coupling so as to have
the same value of $K$, evidently the recursion (\ref{eqn:h_recursion})
preserves this stochastic order, which thus applies for all $t$.
In the limit $t \to\infty$ we thus deduce that
$Q^{*}\preceq Q^{**}\preceq Q^{*}_+$.
Since $\rho $ has finite first moment,
by (\ref{eq:mg-root}) of Theorem \ref{theorem:TreeDecay},
$\mathbb{E}|\tanh(h^{(t)}_+)-\tanh(h^{(t)})|\to0$ as
$t \to\infty$.
Thus, the expectation of the monotone increasing
continuous and bounded function $\tanh(h)$ is the same under
both $Q^{*}$ and $Q^{*}_+$. Necessarily this is also the
expectation of $\tanh(h)$ under $Q^{**}$ and the
uniqueness of the nonnegative fixed point of (\ref{eqn:h_recursion})
follows.
\end{pf*}

We next control the dependence on $\beta$ of the
distribution of the fixed point $h^*$ from Lemma \ref{lemma:Recursive}.
\begin{lemma}\label{lemma:Beta}
Let $\|X-Y\|_{\mathrm{MK}}$ denote the
Monge--Kantorovich--Wasserstein distance
between given laws of random variables $X$ and $Y$
(that is, the infimum of $\mathbb{E}|X-Y|$ over all couplings of $X$
and $Y$).
For any $B>0$ and $\beta_{\max}$ finite there exists a
constant $C=C(\beta_{\max},B)$ such that
if $h^*_{\beta_1}$, $h^*_{\beta_2}$ are the fixed
points of the recursion (\ref{eqn:h_recursion})
for $0 \le\beta_1, \beta_2 \le\beta_{\max}$, then
%
\begin{eqnarray}
\|\tanh(h^*_{\beta_2})-\tanh(h^*_{\beta_1})\|_{\mathrm{MK}}\le C
|\beta
_2-\beta_1|.
\end{eqnarray}
\end{lemma}

\begin{pf} Fixing a random tree $\mathsf{T}= \mathsf{T}(\rho
,\infty)$ of
degree distribution $\rho $, recall that while proving
Lemma \ref{lemma:Recursive} we provided a coupling
of the random variables $\tanh(h^*_\beta)$ and the
Ising root magnetizations $m^{t,+/0}(\beta,B)$ at $\beta$
such that
\[
m^{t,0}(\beta,B) \leq\tanh(h^*_{\beta}) \leq m^{t,+}(\beta,B)
\]
for each $\beta$ and all $t$. By Griffiths inequality
the magnetizations at the root are nondecreasing in $\beta$
so from the bound (\ref{eq:mg-root}) we get
that for $M=M(\beta_{\max},B,\overline{\rho})$ and
any $\beta_1 \le\beta_2 \le\beta_{\max}$,
\begin{eqnarray*}
\mathbb{E}|\tanh( h^*_{\beta_2})-\tanh(h^*_{\beta_1})| &\le&
\mathbb{E}m^{t,0}(\beta_2,B)- \mathbb{E}m^{t,0}(\beta_1,B) +
\frac{M}{t} \\
&\le& (\beta_2-\beta_1) \sup_{\beta\le\beta_{\max}}
\mathbb{E}\biggl\{\frac{\partial m^{t,0}}{\partial\beta}\biggr\}   +
\frac{M}{t}   ,
\end{eqnarray*}
where the expectations
are over the random tree $\mathsf{T}(\rho ,\infty)$.
Considering $t \to\infty$ it thus suffices to
show that $\mathbb{E}[\partial m^{\ell,0}/\partial\beta]$
is bounded, uniformly in $\ell$ and $\beta\le\beta_{\max}$.
To this end, a straightforward calculation yields
\begin{eqnarray*}
\frac{\partial m^{\ell,0}}{\partial\beta}(\beta,B) =
\sum_{(i,j)\in\mathsf{T}(\ell)}
(\langle x_{\o}x_ix_j\rangle- \langle x_{\o}\rangle\langle
x_ix_j\rangle)  ,
\end{eqnarray*}
with $\langle  \cdot  \rangle$ denoting the expectation with
respect to
the Ising measure $\mu^{\ell,0}$. If $i$ is on the path
in $\mathsf{T}(\ell)$ between the root and $j$, then
under the measure $\mu^{\ell,0}$ the variables
$x_{\o}$ and $x_j$ are conditionally independent given $x_i$.
Further, as $x_i \in\{-1,1\}$ it is easy to check that in this case
\[
\langle x_{\o}x_ix_j\rangle- \langle x_{\o}\rangle\langle
x_ix_j\rangle =
\gamma\langle x_{\o};x_i\rangle  ,
\]
where $\gamma$ is the arithmetic
mean of the conditional expected value of $x_j$ for $x_i=-1$
and the conditional expected value of $x_j$ for $x_i=1$.
Thus, $|\gamma| \le1$ and recalling (\ref{eq:ident-corr}) that
$\langle x_{\o};x_i\rangle$ is nonnegative by Griffiths inequality, we
deduce that
\[
\frac{\partial m^{\ell,0}}{\partial\beta}(\beta,B) \le
\sum_{i\in\mathsf{T}(\ell-1)} \Delta_i
\langle x_{\o};x_i\rangle =
\sum_{r=0}^{\ell-1} V_{r,\ell}  ,
\]
where $\Delta_i$ denotes the offspring number at $i \in\mathsf{T}$ and
by (\ref{eq:ident-corr})
\begin{eqnarray*}
V_{r,\ell} \equiv\sum_{i\in\partial\mathsf{T}(r)} \Delta_i
\langle x_{\o};x_i\rangle
= \sum_{i\in\partial\mathsf{T}(r)} \Delta_i.
\partial_{B_i} m^{\ell}(\underline{B},0) |_{\underline{B}=B}
\end{eqnarray*}
[with $m^k(\underline{B},\underline{H})$ the root magnetization for the
measure $\mu^{\underline{B},\underline{H}}$ of (\ref{eq:muBH})].
In view of Lemma \ref{lemma-observation}
we have that $m^k(\underline{B},0)=m^{k-1}(\underline{B},\underline
{H})$ for some
nonnegative vector $\underline{H}$.
By GHS inequality we deduce that for any $i \in\mathsf{T}(k-1)$
\[
\partial_{B_i} m^{k}(\underline{B},0) = \partial_{B_i}
m^{k-1}(\underline{B},\underline{H})
\le\partial_{B_i} m^{k-1}(\underline{B},0)  .
\]
Consequently, $V_{r,\ell}$ is nonincreasing in $\ell$ and
\[
\mathbb{E}\biggl[ \frac{\partial m^{\ell,0}}{\partial\beta} \biggr]
\le\sum_{r=0}^{\ell-1} \mathbb{E}V_{r,\ell}
\leq\sum_{r=0}^{\ell-1} \mathbb{E}V_{r,r}
\leq\sum_{r=0}^{\infty} \mathbb{E}V_{r,r}  .
\]
Further, $m^r(\underline{B},0)$ is independent of the
offspring numbers at $\partial\mathsf{T}(r)$ whose expectation
with respect to the random tree $\mathsf{T}(\rho,\infty)$ is
$\overline{\rho}$. Thus, applying (\ref{eq:ExponentialDecay})
of Corollary \ref{coro:Susc} for $\ell=r$,
$\mathsf{T}=\mathsf{T}(\rho ,\infty)$ and constant magnetic field,
we find that
for some $A$ finite, $\lambda>0$, any $r \ge0$ and all
$\beta\le\beta_{\max}$
\[
\mathbb{E}V_{r,r}
= \overline{\rho}  \mathbb{E}\biggl[ \sum_{i\in\partial\mathsf{T}(r)}
 \partial_{B_i} m^{r}(\underline{B},0) |_{\underline
{B}=B} \biggr]
= \overline{\rho}  \mathbb{E}\biggl[ \sum_{i \in\partial\mathsf{T}(r)}
\langle x_{\o};x_i\rangle \biggr]
\le\overline{\rho}  A e^{-\lambda r}  .
\]
Summing over $r$ gives us the required
uniform boundedness of $\mathbb{E}[\partial m^{\ell,0}/\partial\beta]$
in $\ell$ and $\beta\le\beta_{\max}$.
\end{pf}

%
\section{Algorithms}\label{sec:Algo}

The theorems stated in Section \ref{sec:AlgoImpl} are in fact
consequences of Corollary \ref{coro:Susc}.

\begin{pf*}{Proof of Theorem \ref{theorem:BPConvergence}}
The proof is based on the well-known representation of the iteration
(\ref{eq:BPIteration}) in terms of ``computation tree'' \cite{CompTree}.
Namely, $\nu^{(t)}_{i\to j}( \cdot )$ coincides with the marginal
at the root of the Ising model (\ref{eq:IsingModel}) on a properly
constructed, deterministic tree $\mathsf{T}^{\mathsf{c}}_{i\to j}(t)$ of $t$ generations.
While we refer to the literature for the precise definition of
$\mathsf{T}^{\mathsf{c}}_{i\to j}(t)$, here are some immediate properties:
\begin{longlist}[(a)]
\item[(a)] One can construct an infinite tree $\mathsf{T}^{\mathsf{c}}_{i\to j}(\infty)$
such that, for any $t$, $\mathsf{T}^{\mathsf{c}}_{i\to j}(t)$ is the subtree formed
by the first $t$ generations of $\mathsf{T}^{\mathsf{c}}_{i\to j}(\infty)$.
\item[(b)] The maximal degree of $\mathsf{T}^{\mathsf{c}}_{i\to j}(\infty)$ is
bounded by
the maximal degree of $G$ (and equal to the latter when
$G$ is connected).
\item[(c)] A positive initialization corresponds to adding
$H_{l\to k} = \operatorname{atanh}(\nu^{(0)}_{l\to k}(+1) - \nu
^{(0)}_{l\to k}(-1))$
nonnegative to the field $B$ on the $t$th generation vertices of
$\mathsf{T}^{\mathsf{c}}_{i\to j}(t)$.
\end{longlist}
Denote by $\nu^{+,(t)}_{i\to j}(  \cdot  )$,
$\nu^{0,(t)}_{i\to j}(  \cdot  )$ the messages obtained under
initializations\break
$\nu^{+,(0)}_{k\to l}(+1 )=1$ and
$\nu^{0,(0)}_{k\to l}(+1) =\nu^{0,(0)}_{k\to l}(-1)=1/2$,
respectively.
By Griffiths inequality,
$\nu^{+,(t)}_{i\to j}(+1)$ is nonincreasing in $t$,
$\nu^{0,(t)}_{i\to j}(+1)$ is nondecreasing in $t$
and any positive initialization results with
$\nu^{(t)}_{i \to j}(  \cdot )$ such that
\[
\nu^{+,(t)}_{i\to j}(+1)\ge\nu^{(t)}_{i\to j}(+1)\ge
\nu^{0,(t)}_{i\to j}(+1)  .
\]
By Corollary \ref{coro:Susc} we have that
$\nu^{+,(t)}_{i\to j}(+1)-\nu^{0,(t)}_{i\to j}(+1) \le A
e^{-\lambda t}$
for all $t \ge0$. Since $A<\infty$ and $\lambda>0$ depend
only on $\beta$, $B$ and the maximal degree of $G$, this
immediately yields our thesis.
\end{pf*}

\begin{pf*}{Proof of Theorem \ref{theorem:BPCorrectness}}
We use an additional property of the computation tree:
\begin{longlist}
\item[(d)] If $\mathsf{B}_{i}(k)$ is a tree then $\mathsf{T}^{\mathsf{c}}_{i
\to j}(k)$
is a tree rooted at $i \to j$ whose vertices
are the directed edges on the maximal subtree of $\mathsf{B}_{i}(k)$
rooted at $i$ that does not include~$j$.
\end{longlist}
Without loss of generality we may and shall assume that
$t>r$. For $U=\mathsf{B}_{i_*}(r)$ consider the local marginal approximations
$\nu_U^{+}( \cdot )$, $\nu_U^{0}( \cdot )$ defined as in
(\ref{eq:LocalMarg}) except that the fixed point messages
$\nu^*_{i\to j(i)}( \cdot )$ at $i \in\partial\mathsf{B}_{i_*}(r)$
are replaced by\vspace{-2pt} those obtained after $(t-r)$ iterations starting at
$\nu^{+,(0)}_{k\to l}(+1 )=1$ and
$\nu^{0,(0)}_{k\to l}(+1) =\nu^{0,(0)}_{k\to l}(-1)=1/2$,
respectively. Since $\mathsf{B}_{i_*}(t)$ is a tree, here
$j(i)$ is necessarily the neighbor of $i$
on the path from $i_*$ to $i \in\partial\mathsf{B}_{i_*}(r)$ and
from the preceding property (d) we see that
$\mathsf{T}^{\mathsf{c}}_{i \to j(i)}(t-r)$ corresponds to the subtree of $i$ and its
lines of descendant in $\mathsf{B}_{i_*}(t)$. By property (c)
we thus have that $\nu_U^{+} ( \cdot )$ and
$\nu_U^{0} ( \cdot )$ are
the marginals on $U$ of the Ising model $\nu^+$ on $G$
with $B_i = \infty$ at all $i \notin\mathsf{B}_{i_*}(t)$ and
the Ising model $\nu^0$ on the vertices of $G$
and the edges within the tree $\mathsf{B}_{i_*}(t)$.
Such reasoning also shows that the probability measure
$\nu_U$ of (\ref{eq:LocalMarg}) is the marginal on $U$
of the Ising model $\nu$ on vertices of $G$
and edges of $\mathsf{B}_{i_*}(t)$ with an additional
nonnegative magnetic field
$H_{l\to k} = \operatorname{atanh}(\nu^*_{l\to k}(+1) - \nu^*_{l\to k}(-1))$
at $\partial\mathsf{B}_{i_*}(t)$. Consequently,
with $x_F \equiv\prod_{i \in F} x_i$
we have by Griffiths inequality that for any $F \subseteq U$
\[
\langle \nu^0, x_F \rangle \le\langle \nu, x_F \rangle \le
\langle \nu^+ , x_F \rangle  ,
\qquad
\langle \nu^0, x_F \rangle \le\langle \mu, x_F \rangle \le
\langle \nu^+ , x_F \rangle  ,
\]
and we deduce that for any $F \subseteq U$,
\[
|\langle \mu, x_F \rangle - \langle \nu, x_F \rangle| \le
\langle \nu^+, x_F \rangle - \langle \nu^0, x_F \rangle
\le2 \|\nu^+_U - \nu^0_U\|_{\mathrm{TV}}  .
\]
Recall that since $x_i \in\{-1, 1 \}$,
for any possible value $\underline{y}=\{ y_i, i \in U \}$ of
$\underline{x}_U$,
\[
\mathbb{I}(\underline{x}_U = \underline{y}) = 2^{-|U|} \prod_{i \in
U} (1+ y_i x_i)
= 2^{-|U|} \sum_{F \subseteq U} y_F   x_F  ,
\]
and with $|y_F| \le1$ it follows that
\begin{eqnarray*}
|\mu_U(\underline{y}) - \nu_U(\underline{y})| &=&
2^{-|U|} \biggl| \sum_{F \subseteq U} y_F (\langle \mu_U, x_F \rangle -
\langle \nu_U,
x_F \rangle) \biggr|
\\
&\le& \max_{F \subseteq U} |\langle \mu_U, x_F \rangle - \langle \nu
_U, x_F \rangle|
\le2 \|\nu^+_U - \nu^0_U\|_{\mathrm{TV}}  .
\end{eqnarray*}
This applies for any of the $2^{|U|}$ possible values of $\underline
{x}_U$, so
\begin{eqnarray*}
\|\mu_U( \cdot ) - \nu_U( \cdot )\|_{\mathrm{TV}} \le
2^{|U|}   \|\nu^+_U( \cdot ) - \nu^0_U( \cdot )\|_{\mathrm{TV}}  .
\end{eqnarray*}
Applying Corollary \ref{coro:Susc} for the
deterministic tree $\mathsf{B}_{i_*}(t)$ rooted
at $i_*$, we get the bound (\ref{eq:bd-influence})
on the right side of the preceding inequality with
$\delta(k)=A \exp(-\lambda k)$, some finite $A$ and
$\lambda>0$ that depend only on $\beta$, $B$ and $\Delta$.
Thus, noting that $|U|= |\mathsf{B}_{i_*}(r)| \le\Delta^{r+1} + 1$
we establish our thesis upon choosing $c=c(A,C,\Delta)$
large enough.
\end{pf*}
%
%
\section{From trees to graphs}\label{sec:Graphs}

We start with the following technical lemma.
\begin{lemma}\label{lemma:expansion}
Consider a convex set $\mathcal{K}\subseteq{\mathbb{R}}$ and symmetric
twice differentiable functions $F_\ell\dvtx \mathcal{K}^\ell\to{\mathbb{R}}$
with $F_0$ constant, such that for some finite constant~$c$,
\[
\sup_{\ell} \sup_{\mathcal{K}^\ell} \biggl|
\frac{\partial^2 F_\ell}{\partial x_1\, \partial x_2}\biggr| \le2 c  .
\]
Suppose i.i.d. $X, X_i \in\mathcal{K}$
are such that $\ell^{-1} \mathbb{E}|\partial_{x_1} F_\ell
(x,X_2,\ldots
,X_\ell)|$
is bounded uniformly in $\ell$ and $x \in\mathcal{K}$ and the
independent, square-integrable,
nonnegative integer valued random variable $L$ satisfies
%
\begin{eqnarray}\label{eq:vanish-first-der}
\mathbb{E}[ L \partial_{x_1} F_L (x,X_2,\ldots,X_L) ] = 0 \qquad
\forall x\in\mathcal{K} .
\end{eqnarray}
Then, for any i.i.d. $Y,Y_i \in\mathcal{K}$ also independent of $L$,
\begin{eqnarray}\label{eq:Claim}
&&|\mathbb{E}[ F_L(Y_1,\dots,Y_L)- F_L(X_1,\dots,X_L) ]|\nonumber
\\[-8pt]\\[-8pt]
&&\qquad\le c
\mathbb{E}[L
(L-1)]
\|X-Y\|^2_{\mathrm{MK}}  .\nonumber
\end{eqnarray}
\end{lemma}

\begin{pf} Our thesis trivially holds if either
$\|X-Y\|_{\mathrm{MK}}=0$ or $\|X-Y\|_{\mathrm{MK}}=\infty$, so
without loss of generality, fixing $\gamma>1$ we assume hereafter that
$(X_i,Y_i)$ are i.i.d. pairs, independent on $L$ and
coupled in such a way that $\mathbb{E}|X_i-Y_i|\le\gamma
\|X-Y\|_{\mathrm{MK}}$ is
finite. It is easy to check that almost surely,
\begin{eqnarray}\label{eq:taylor}
&&F_\ell(Y_1,\dots,Y_\ell)-F_\ell(X_1,\dots,X_\ell)\nonumber
\\[-8pt]\\[-8pt]
&&\qquad =\sum
_{i=1}^\ell\Delta_i
F_\ell+ \sum_{i\neq j}^\ell f_{ij}^{(\ell)} (Y_i-X_i) (Y_j-X_j)  ,\nonumber
\end{eqnarray}
where $\Delta_i F_\ell= (Y_i-X_i) \int_0^1 \partial_{x_i}
F_\ell(X_1,\dots,t Y_i + (1-t) X_i,\dots, X_\ell)\, \mathrm{d}t$ and
each of the terms
\begin{eqnarray*}
f_{ij}^{(\ell)} &=& \int_{0}^1 \int_{0}^t
\frac{\partial^2 F_\ell}{\partial x_i \,\partial x_j}
\bigl(s Y_1 +(1-s)X_1,\ldots,
\\
&&\hspace{70pt}{}tY_i+(1-t)X_i,\ldots,sY_\ell+(1-s)X_\ell\bigr)
\,\mathrm{d}
s\, \mathrm{d}t  ,
\end{eqnarray*}
is bounded by $c$. For i.i.d. $(X_i,Y_i)$, by the
symmetry of the functions $F_\ell$ with respect to their arguments,
the assumed boundedness of
$\ell^{-1} \mathbb{E}|\partial_{x_1} F_\ell(x,X_2,\break\ldots,X_\ell)|$
implies integrability of $\Delta_i F_\ell$ with
$\mathbb{E}\Delta_i F_\ell$ independent of $i$ and
$\ell^{-1} \mathbb{E}|\Delta_i F_\ell|$ uniformly bounded. This in
turn implies the integrability of $\sum_{i=1}^L \Delta_iF_L$
for any $L$ square integrable and independent of $(X_i,Y_i)$, so
by Fubini's theorem and our assumption (\ref{eq:vanish-first-der}),
\begin{eqnarray*}
&&\mathbb{E}\Biggl[ \sum_{i=1}^L \Delta_i F_L \Biggr]
\\
&&\qquad = \mathbb{E}
[L \Delta_1
F_L]
\\
&&\qquad= \mathbb{E}\biggl[ (Y_1-X_1)
\\
&&{}\qquad\qquad \times\int_0^1\mathbb{E}\bigl[L \partial_{x_1} F_L
 \bigl(tY_1+(1-t)X_1,X_2, \dots, X_L\bigr)
  |
X_1,Y_1 \bigr]\, \mathrm{d}t \biggr] = 0  .
\end{eqnarray*}
Thus, considering the expectation of (\ref{eq:taylor}), by
the uniform boundedness of $f_{ij}^{(\ell)}$ and the
independence of $L$ on the i.i.d. pairs $(X_i,Y_i)$,
we deduce that
\begin{eqnarray*}
| \mathbb{E}[ F_L (Y_1,\dots,Y_L)- F_L (X_1,\dots,X_L)
] |
&\le& c \mathbb{E}\sum_{i\neq j}^L |Y_i-X_i| |Y_j-X_j|
\\
&\le&
\gamma^2 c \mathbb{E}[L (L-1)]   \|X-Y\|^2_{\mathrm{MK}}  .
\end{eqnarray*}
Finally, taking $\gamma\downarrow1$ yields the bound (\ref{eq:Claim}).
\end{pf}

\begin{remark}\label{rem-mod} It is not hard to adapt the proof of
the lemma so as to replace $F_1\dvtx \mathcal{K}\mapsto{\mathbb{R}}$ by
$0.5 F_1(x,y)$ for a twice differentiable
symmetric function $F_1 \dvtx \mathcal{K}^2 \mapsto{\mathbb{R}}$.
Taking $P_\ell=\mathbb{P}(L=\ell)$ the contribution of $L=1$ to
the left-hand side of (\ref{eq:vanish-first-der}) is then
$P_1 \mathbb{E}[ \partial_{x_1} F_1 (x,X_2) ]$ and
the bound (\ref{eq:Claim}) is modified to
%
\begin{eqnarray}\label{eq:Claim-mod}
&&\biggl|\frac{P_1}{2} \mathbb{E}[ F_1(Y_1,Y_2) - F_1(X_1,X_2)]\nonumber
\\
&&{}\quad +
\sum_{\ell\ge2} P_\ell
\mathbb{E}[ F_\ell(Y_1,\dots,Y_\ell)- F_\ell(X_1,\dots,X_\ell)
]\biggr|
\\
&&\qquad\le c
\mathbb{E}[L^2]
  \|X-Y\|^2_{\mathrm{MK}}  .\nonumber
\end{eqnarray}
\end{remark}

Consider the functional
$h\mapsto\varphi_h$ that, given a random variable $h$,
evaluates the right-hand side of Equation~(\ref{eqn:phi}).
It is not hard to check that $\varphi_h$
is well defined and finite for every random variable $h$.
The following corollary of Lemma \ref{lemma:expansion}
plays an important role in the proof of Theorem \ref{theorem:free_energy}.
\begin{coro}\label{lemma:Stationarity}
There exist nondecreasing
finite $c(|\beta|)$ such that if $\overline{\rho}<\infty$ and
$h^*$ is a fixed point of the distributional identity
(\ref{eqn:h_recursion}) for some $\beta,B \in{\mathbb{R}}$ then
%
\begin{eqnarray}
|\varphi_h(\beta,B)-\varphi_{h^*}(\beta,B)| \le c(|\beta|)
\overline{P}  \overline{\rho}
\|\tanh(h)-\tanh(h^*)\|^2_{\mathrm{MK}}  .
\end{eqnarray}
%
\end{coro}

\begin{pf} Setting $u=\tanh(\beta)$ so $|u|<1$, we
verify the conditions of\break Lemma~\ref{lemma:expansion} when
$X_i$ are i.i.d. copies of
$X=\tanh(h^*)$ and $Y_i$ i.i.d. copies of $Y = \tanh(h)$, all of
whom take values in $\mathcal{K}=[-1,1]$ and
are independent of the random variable $L$. We apply the lemma
in this setting for the symmetric, twice differentiable functions
\begin{eqnarray*}
F_\ell(x_1,\dots,x_\ell) &=& -\frac{1}{(\ell-1)}
\sum_{1 \le i<j \le\ell} \log(1+u x_i x_j)
\\
&&{}+\log\Biggl\{ e^B \prod_{i=1}^{\ell}(1 + u x_i)
+ e^{-B} \prod_{i=1}^{\ell} (1 - u x_i)\Biggr\}
\end{eqnarray*}
for $\ell\ge2$, and as in Remark \ref{rem-mod},
\begin{eqnarray*}
F_1(x_1,x_2) &=& - \log(1+u x_1 x_2)
+\log\{ e^B (1 + u x_1) + e^{-B} (1 - u x_1)\}
\\
&&{}+\log\{ e^B (1 + u x_2) + e^{-B} (1 - u x_2)\}  .
\end{eqnarray*}
Indeed, setting $\psi(x,y)=u y/(1+u x y)$ and for each $\ell\ge1$
%
\begin{equation}\label{eq:gdef}
g_\ell(x_2,\ldots,x_\ell)=\tanh\Biggl(B+\sum_{j=2}^\ell\operatorname{atanh}(u
x_j)\Biggr)
\end{equation}
[so $g_1=\tanh(B)$], it is not hard to verify that
$\partial_{x_1} F_1 (x_1,x_2) = \psi(x_1,g_1) - \psi(x_1,x_2)$
while for $\ell\ge2$
%
\begin{equation}\label{eq:first-der}
\partial_{x_1} F_\ell(x_1,\dots,x_\ell)
= \psi(x_1,g_\ell(x_2,\dots,x_\ell))
- \frac{1}{\ell-1} \sum_{j=2}^\ell\psi(x_1,x_j)  .
\end{equation}
In particular, $g_\ell( \cdot )$ are
differentiable functions from $\mathcal{K}^{\ell-1}$ to $\mathcal
{K}$, such that
$\partial_{x_2} g_\ell$
are uniformly bounded [by $a = |u|/(1-u^2)$] and
$\partial_y \psi(x,y)$ is uniformly bounded
on $\mathcal{K}^2$ [by $b = |u|/(1-|u|)^2$]. Consequently,
$\partial_{x_1} F_\ell$ and
$\partial^2 F_\ell/\partial x_1\, \partial x_2$ are also uniformly
bounded [by $2/(1-|u|)$ and $b(a+1)=2c(|\beta|)$, respectively].
Further, $h^*$ is a fixed point of (\ref{eqn:h_recursion}), hence
$X_1 \stackrel{\mathrm{d}}{=}g_K (X_2,\dots,X_K)$.
With $X_i$ identically distributed and $\overline{P}\rho_k = k P_k$
we thus find as required in (\ref{eq:vanish-first-der}) that
\begin{eqnarray}\label{eq:fixpt}
&&P_1 \mathbb{E}[ \partial_{x_1} F_1 (x,X_2) ]  + \sum
_{k \ge2}
k P_k \mathbb{E}[ \partial_{x_1} F_k (x,X_2,\ldots,X_k) ] \nonumber
\\[-8pt]\\[-8pt]
&&\qquad =
\overline{P}  \Biggl\{   \sum_{k=1}^\infty\rho_k \mathbb{E}[ \psi
(x,g_k(X_2,\dots,X_k)) ]
- \mathbb{E}\psi(x,X_1)   \Biggr\} = 0  .\nonumber
\end{eqnarray}
Noting that $\mathbb{E}[L^2] = \overline{P}  \overline{\rho} $
our thesis is merely
the bound (\ref{eq:Claim-mod}) upon confirming that
\begin{eqnarray*}
\varphi_h &=& F_0 + \frac{P_1}{2} \mathbb{E}F_1(Y_1,Y_2) +
\sum_{\ell\ge2} P_\ell  \mathbb{E}F_\ell(Y_1,\dots,Y_\ell),
\\
\varphi_{h^*} &=&
F_0 + \frac{P_1}{2} \mathbb{E}F_1(X_1,X_2) +
\sum_{\ell\ge2} P_\ell  \mathbb{E}F_\ell(X_1,\dots,X_\ell)
\end{eqnarray*}
for some constant $F_0$
and that both series are absolutely summable.
\end{pf}

Let $\overline{\mathsf{T}}(\rho ,\infty)$ denote
the infinite random tree obtained by ``gluing''
two independent trees from the ensemble $\mathsf{T}(\rho ,\infty)$ through
an extra edge $e$ between their roots and considering~$e$ as the root of
$\overline{\mathsf{T}}(\rho ,\infty)$ denote by $\overline{\mathsf
{T}}(\rho ,t)$ the subtree
formed by its first $t$ generations [i.e., consisting of $e$ and the
corresponding two independent copies from $\mathsf{T}(\rho ,t)$].
An alternative way to sample from $\overline{\mathsf{T}}(\rho
,\infty)$ is to have
independent offspring number $k-1$ with probability
$\rho _k$ at each end of the root edge~$e$ and thereafter
independently sample from this offspring distribution
at each revealed new node of the tree.
Equipped with these notations we have the following consequence of
the local convergence of the graph sequence $\{G_n\}$.
\begin{lemma}\label{lemma:EdgeTree}
Suppose a uniformly sparse graph sequence $\{G_n\}$
converges locally to the random tree $\mathsf{T}(P,\rho ,\infty)$.
Fixing a nonnegative integer $t$, for each $(i,j)\in E_n$ denote
the subgraph of $G_n$ induced by vertices at distance at
most $t$ from $(i,j)$ by $\mathsf{B}_{ij}(t)$. Let $F(\cdot)$ be a
fixed, bounded function on the collection
of all possible subgraphs that may occur as $\mathsf{B}_{ij}(t)$,
such that $F(T_1)=F(T_2)$ whenever $T_1 \simeq T_2$.
Then,
%
\begin{eqnarray}\label{eq:edge-conv}
\lim_{n\to\infty}\frac{1}{n}\sum_{(i,j)\in E_n} F(\mathsf
{B}_{ij}(t)) =
\frac{\overline{P}}{2}  \mathbb{E}\{ F(\overline{\mathsf
{T}}(\rho ,t))\}  .
\end{eqnarray}
\end{lemma}

\begin{pf} Denoting by $\mathbb{E}_{(ij)}(\cdot)$ the expectation
with respect to a\break uniformly chosen edge $(i,j)$ in $E_n$,
the left-hand side of (\ref{eq:edge-conv}) is merely
$(|E_n|/n)\mathbb{E}_{(ij)}\{F(\mathsf{B}_{ij}(t))\}$.
A~uniformly chosen edge can be sampled by first selecting
a vertex $i$ with probability proportional to its degree
$|{\partial i}|$ and then picking one of its neighbors $j=j(i)$ uniformly.
Thus, denoting by $\mathbb{E}_n(\cdot)$ the expectation with respect to
a uniformly chosen random vertex $i \in[n]$, we have that
\begin{eqnarray*}
\mathbb{E}_{(ij)}\{F(\mathsf{B}_{ij}(t))\} = \frac{\mathbb{E}_n\{
|{\partial i}|  F(\mathsf{B}
_{ij(i)}(t))\}}
{\mathbb{E}_n\{|{\partial i}|\}}  .
\end{eqnarray*}
Marking uniformly at random one offspring of $\o$
in $\mathsf{T}(P,\rho ,t+1)$ [as corresponding to $j(i)$],
let $\mathsf{T}_* (t+1)$ denote the subtree
induced by vertices whose distance from either
$\o$ or its marked offspring is at most $t$.
Since $\mathsf{B}_{ij(i)}(t) \subseteq\mathsf{B}_i(t+1)$ and
with probability $q_{t,k} \to1$ as $k \to\infty$
the random tree $\mathsf{T}(P,\rho ,t+1)$
belongs to the finite collection of trees with
$t+1$ generations and maximal degree at most $k$,
it follows by dominated convergence and the
local convergence of $\{G_n\}$ that for any fixed $l$,
\begin{eqnarray*}
&&\lim_{n \to\infty}
\mathbb{E}_n\bigl[  |{\partial i}|  \mathbb{I}(|{\partial i}|\le
l) F\bigl(\mathsf{B}_{ij(i)}(t)\bigr) \bigr]
\\
&&\qquad = \mathbb{E}_\rho \bigl\{\Delta_{\o} \mathbb{I}(\Delta_{\o} \le l)
F\bigl(\mathsf{T}_* (t+1)\bigr)\bigr\}  ,
\end{eqnarray*}
where $\mathbb{E}_\rho (\cdot)$ and
$\Delta_{\o}$ denote expectations and
the degree of the root, respectively, in
$\mathsf{T}(P,\rho ,\infty)$. Similarly,
\[
\lim_{n \to\infty}
\mathbb{E}_n\{  |{\partial i}|   \mathbb{I}(|{\partial i}| \le l)
  \}
= \mathbb{E}_\rho \Delta_{\o} \mathbb{I}(\Delta_{\o} \le l)  .
\]
Further, by the uniform sparsity of $\{G_n\}$,
%
\begin{eqnarray*}
&&\limsup_{n \to\infty}
\bigl|\mathbb{E}_n\bigl[  |{\partial i}|  \mathbb{I}(|{\partial i}|
> l) F\bigl(\mathsf{B}_{ij(i)}(t)\bigr)
\bigr] \bigr|
\\
&&\qquad \le\|F\|_\infty
\limsup_{n \to\infty}
\mathbb{E}_n[  |{\partial i}|  \mathbb{I}(|{\partial i}| > l)
]
\end{eqnarray*}
goes to zero as $l \to\infty$.
Since $P$ has a finite first moment,
$\Delta_{\o}$ is integrable, so by the preceding, upon
taking $l \to\infty$ we deduce by dominated convergence that
\begin{eqnarray*}
\lim_{n\to\infty}\mathbb{E}_{(ij)}\{F(\mathsf{B}_{ij}(t))\} =
\frac{
\mathbb{E}_\rho \{\Delta_{\o}F(\mathsf{T}_*(t+1))\}}{\mathbb
{E}_\rho
\{\Delta_{\o}\}}  .
\end{eqnarray*}
To complete the proof note
that the right-hand side of the last expression is precisely
$\mathbb{E}\{ F(\overline{\mathsf{T}}(\rho ,t))\}$ and we have also
shown that
$2 |E_n|/n = \mathbb{E}_n \{|{\partial i}|\} \to\mathbb{E}_\rho
\Delta_{\o} = \overline{P}$.
\end{pf}

\begin{pf*}{Proof of Theorem \ref{theorem:free_energy}}
Since $\phi_n(\beta,B) \equiv\frac{1}{n}\log Z_n(\beta,B)$
is invariant under $B\to-B$ and is uniformly (in $n$)
Lipschitz continuous in $B$ with Lipschitz constant one,
it suffices to fix $B>0$ and show that $\phi_n(\beta,B)$
converges as $n \to\infty$ to
the predicted $\varphi_{h^*}(\beta,B)$ of (\ref{eqn:phi}),
whereby $h^*=h^*_\beta$ is the unique fixed point of
the recursion (\ref{eqn:h_recursion}) that is supported on
$[0,\infty)$ (see Lemma \ref{lemma:Recursive}).

This is obviously true for $\beta=0$ since
$\phi_n(0,B) = \log(2\cosh B) = \varphi_{h} (0,B)$.
Next, denoting by $\langle   \cdot  \rangle_n$ the expectation
with respect to the Ising measure on $G_n$ (at
parameters $\beta$ and $B$), it is easy to see that
%
\begin{eqnarray}\label{eq:DerivativeGraph}
\partial_\beta\phi_n(\beta,B) = \frac{1}{n}\sum_{(i,j)\in E_n}
\langle x_i x_j\rangle_n .
\end{eqnarray}
Clearly $|\partial_\beta\phi_n(\beta,B)| \le|E_n|/n$
is bounded by the uniform sparsity of $\{G_n\}$ so
it is enough to show that
the expression in (\ref{eq:DerivativeGraph})
converges to the partial derivative of
$\varphi_{h^*_\beta} (\beta,B)$ with respect to $\beta$.
Turning to compute the latter derivative, by
Lemma~\ref{lemma:Beta} and Corollary \ref{lemma:Stationarity}
we can ignore the dependence of $h^*_\beta$ on $\beta$.
That is, we simply compute the partial derivative in $\beta$ of the
expression (\ref{eqn:phi}) while
considering (the law of) $h_i$ to be fixed. Indeed,
with notation $u=\tanh(\beta)$ and $X_i = \tanh(h_i)$
as in the derivation of Corollary \ref{lemma:Stationarity},
a direct computation leads by the exchangeability of $X_i$ to
\begin{eqnarray*}
\partial_\beta  \varphi(\beta,B) &=&
\frac{\overline{P}}{2} u -
\frac{\overline{P}}{2} (1-u^2) \mathbb{E}[ \psi(X_1,X_2) ]
\\
&&{}+ (1-u^2) \mathbb{E}[ L \psi(X_1,g_L(X_2,\ldots,X_L))]
\end{eqnarray*}
for $\psi(x,y)=xy/(1+uxy)$ and $g_\ell(x_2,\ldots,x_\ell)$ of
(\ref{eq:gdef}). Further, the
fixed point property (\ref{eq:fixpt})
applies for any bounded measurable $\psi(\cdot)$, so we deduce
that
\begin{eqnarray*}
\mathbb{E}[L \psi(X_1,g_L(X_2,\ldots,X_L))] &=&
\overline{P}  \mathbb{E}[\psi(X_1,g_K(X_2,\ldots,X_K))]
\\
&=& \overline{P}  \mathbb{E}[ \psi(X_1,X_2) ]  .
\end{eqnarray*}
Consequently, it is not hard to verify that
%
\begin{eqnarray}\label{eq:FreeDer}
\partial_\beta \varphi(\beta,B) =
\frac{\overline{P}}{2} \mathbb{E}\biggl\{\frac{u + X_1 X_2}{1+ u
X_1 X_2}
\biggr\}
= \frac{\overline{P}}{2} \mathbb{E}[ \langle x_i x_j\rangle
_{\overline{\mathsf{T}}} ]   ,
\end{eqnarray}
where $\langle \cdot\rangle_{\overline{\mathsf{T}}}$ denotes the
expectation with
respect to the Ising model
\[
\mu_{\overline{\mathsf{T}}}(x_i,x_j) = \frac{1}{z_{ij}} \exp
\{
\beta x_ix_j+H_i x_i+H_j x_j\}
\]
on one edge $(ij)$ and random magnetic fields
$H_i$ and $H_j$ that are independent copies of $h^*_\beta$.

In comparison, fixing a positive integer $t$,
by Griffiths inequality the correlation
$\langle x_i x_j\rangle_n$ lies between the correlations
$F_0(\mathsf{B}_{ij}(t)) \equiv\langle x_i x_j\rangle^0_{\mathsf
{B}_{ij}(t)}$
and $F_+(\mathsf{B}_{ij}(t)) \equiv\langle x_i x_j\rangle^+_{\mathsf
{B}_{ij}(t)}$
for the Ising model on the subgraph
$\mathsf{B}_{ij}(t)$ with free and plus, respectively,
boundary conditions at $\partial\mathsf{B}_{ij}(t)$.
Thus, in view of (\ref{eq:DerivativeGraph})
\[
\frac{1}{n}\sum_{(i,j)\in E_n} F_0(\mathsf{B}_{ij}(t)) \le
\partial_\beta  \phi_n(\beta,B) \le
\frac{1}{n}\sum_{(i,j)\in E_n} F_+(\mathsf{B}_{ij}(t))  ,
\]
and taking
$n\to\infty$ we get by Lemma \ref{lemma:EdgeTree} that
\begin{eqnarray*}
\frac{\overline{P}}{2} \mathbb{E}[ F_0(\overline{\mathsf{T}}(\rho
,t)) ] &\le&
\liminf_{n\to\infty} \partial_\beta  \phi_n(\beta,B)
\\
&\le&\limsup_{n\to\infty} \partial_\beta  \phi_n(\beta,B)
\le
\frac{\overline{P}}{2}
\mathbb{E}[ F_+ (\overline{\mathsf{T}}(\rho ,t))]  .
\end{eqnarray*}
To compute $F_{0/+} (\overline{\mathsf{T}}(\rho ,t))$ we first sum over
the values of $x_k$ for
$k \in\overline{\mathsf{T}}(\rho ,t) \setminus\{i,j\}$. This has the
effect of reducing $F_{0/+} (\overline{\mathsf{T}}(\rho ,t))$ to a
form of $\langle x_i x_j\rangle_{\overline{\mathsf{T}}}$. Further, as
shown in the proof of Lemma \ref{lemma:Recursive},
we get $F_{0/+} (\overline{\mathsf{T}}(\rho ,t))$
by setting for $H_i$ and $H_j$ two
independent copies of the variables
$h^{(t)}$ and
$h_{+}^{(t)}$, respectively, which converge in law to $h_\beta^*$
when $t \to\infty$. We also saw there that
the functional $\Psi_U(\nu) = \mathbb{E}[ \langle x_i x_j\rangle
_{\overline{\mathsf{T}}}]$
[for continuous and bounded
$U(H_i,H_j)=(u + \tanh(H_i) \tanh(H_j))/(1+ u \tanh(H_i) \tanh(H_j))$],
is continuous with respect to the weak convergence of
the law $\nu$ of $H_i$. Consequently, by (\ref{eq:FreeDer})
\[
\lim_{t \to\infty} \frac{\overline{P}}{2} \mathbb{E}[ F_{0/+}
(\overline{\mathsf{T}}(\rho
,t)) ] =
\partial_\beta \varphi(\beta,B)  ,
\]
which completes the proof of the theorem.
\end{pf*}
%

%

\printaddresses

\end{document}